\documentclass[a4paper,12pt]{article}
\usepackage[T1]{fontenc}
\usepackage[utf8]{inputenc}
\input{epsf.sty}
\usepackage{hyperref}
\usepackage[final]{epsfig}
\usepackage{color}
\usepackage{cite}
\usepackage{graphicx}
\usepackage{subcaption}
\usepackage{caption}
\usepackage{tikz}
\usetikzlibrary{decorations.pathreplacing}
\usepackage{amsmath, amsthm, amssymb, dsfont, mathrsfs}
\usepackage{comment}
\usetikzlibrary{patterns}

\newtheorem {thm}{Theorem}[section]

\newtheorem {lem}[thm]{Lemma}

\theoremstyle{plain}

\def\gi{\overline\gamma}

\def\N{{\Bbb N}}

\def\Z{{\Bbb Z}}

\def\Q{{\Bbb Q}}

\def\sm{\setminus}
\def\one{{\mathds 1}}

\newcommand{\e}{\mathrm e}
\renewcommand{\d}{\mathrm d}
\newcommand\restr[2]{{
  \left.\kern-\nulldelimiterspace 
  #1 
  \vphantom{\big|} 
  \right|_{#2} 
  }}

\def\0{{\bf 0}}

\def\eps{\varepsilon}

\newcommand{\zd}{{{\mathbb Z}^d}}

\def\eps{\varepsilon}

\def\phi{\varphi}

\def\g{\gamma}

\def\s{\sigma}

\def\o{\omega}

\def\D{\Delta}

\def\L{\Lambda}

\def\O{\Omega}

\def\T{\T}

\def\ot{\tilde{\omega}}

\def\FF{{\cal F}}

\begin{document}

\title{Gibbsianness and non-Gibbsianness for Bernoulli lattice fields under removal of isolated sites}

\author{
Benedikt Jahnel\footnote{
Weierstrass Institute for Applied Analysis and Stochastics, Mohrenstrasse 39
10117 Berlin, Germany
\texttt{Benedikt.Jahnel@wias-berlin.de}}
\, and
Christof K\" ulske\footnote{ Ruhr-Universit\"at Bochum, Fakult\"at f\"ur Mathematik, 44801 Bochum, Germany 
\texttt{Christof.Kuelske@ruhr-uni-bochum.de}}
}

\newcommand{\CC}[1]{{\color{blue} #1}}
\newcommand{\DE}[1]{{\color{red} #1}}
\newcommand{\CK}[1]{{\color{green} #1}}

\maketitle

\begin{abstract}
We consider the i.i.d.~Bernoulli field $\mu_p$ on $\Z^d$ with occupation density $p\in [0,1]$. To each realization of the set of occupied sites we apply a thinning map that removes all occupied sites that are isolated in graph distance. We show that, while this map seems non-invasive for large $p$, as it changes only a small fraction $p(1-p)^{2d}$ of sites, there is $p(d)<1$ such that for all $p\in(p(d),1)$ the resulting measure is a non-Gibbsian measure, i.e., it does not possess a continuous version of its finite-volume conditional probabilities. On the other hand, for small $p$, the Gibbs property is preserved. 
\end{abstract}

\smallskip
\noindent {\bf AMS 2000 subject classification:} primary: 60D05, 60K35; secondary: 82B20
\bigskip 

{\em Keywords: Gibbsianness, Gibbs-uniqueness, Bernoulli field, local thinning, two-layer representation, Dobrushin uniqueness, Peierls' argument} 

\section{Introduction}\label{Sec_Setting}
Random fields under local maps are defined and analyzed in different fields of probability and statistics.  In these studies a random field is typically a stochastic process whose variables take values in a subset of the real numbers, and which has a geometric index set such as the infinite lattice, a graph with infinite vertex set, or the Euclidean space. It is of interest to understand the behavior of the given random field with distribution $\mu$ under application of a map $T$ that acts on infinite-volume realizations of the process, and investigate resulting properties of the image process whose distribution we denote by $\mu'$. 
Relevant deterministic maps $T$ typically have local-dependence windows. They may also be generalized to stochastic kernels that act locally.  In this general setup let us call $\mu$ the {\em first-layer measure}, and the measure $\mu'$, which we shall be mostly interested in, the {\em second-layer measure}. 

This setup is of theoretical interest,  but  also occurs in many applications of natural sciences, engineering, and statistics.  We mention thinning transformations of Poisson point processes (in which points from a realization are omitted according to local rules)~\cite{bremaud1979optimal,brown1979position,isham1980dependent,rolski1991stochastic,last1993dependent,ball2005poisson,moller2010thinning,blaszczyszyn2019determinantal}, transformations in image analysis~\cite{MR1100283,MR1344683}, and renormalization group transformations in statistical mechanics (where a physical system like a ferromagnet is considered on increasingly large scales by means of maps which forget details on small scales)~\cite{griffiths1978position,griffiths1979mathematical,EnFeSo93}.  Processes appearing as local maps in this way can also be viewed as generalizations of the much used {\em hidden Markov models}~\cite{MR1323178}. 

Hidden Markov models  are  images under local kernels of an underlying {\em first-layer} Markov chain, and appear when noisy observations shall be modeled, and so the generalization is made to a situation with spatial index sets, and more complex dependence structures. 

It was discovered first in the context of such renormalization group transformations that strictly local transformations acting on a spatially Markovian random field $\mu$ indexed by a lattice, may result in singularities in the image measure $\mu'$. Two concrete examples for this are provided by the low-temperature Ising model under a block average transformation, or the projection to a sublattice~\cite{MR1012855,EnFeSo93,bricmont1998renormalization,le2013almost,MR3830302}. The original motivation of such renormalization group transformations was, suggested by heuristic schemes of theoretical physics, to understand the iterated {\em coarse-graining dynamics} of the level of Hamiltonians to investigate critical behavior~\cite{wilson1974renormalization}. 

The singularity in the second-layer measure $\mu'$ means, that $\mu'$ loses not only the spatial Markov property of $\mu$ under the map $T$ (which is less surprising), but it even loses the more general {\em Gibbs property}. This is more severe, it means that the infinite second-layer system acquires internal long-range dependence, and in particular does not posses a well-behaved Hamiltonian with good summability properties of its interaction potentials anymore. The singular long-range dependence appears on the level of finite-volume conditional probabilities of the image measure, which are not quasilocal functions of their conditioning. Put equivalently, finite-volume sub-systems depend on their boundary conditions arbitrarily far away, and their behavior cannot be described by kernels that are continuous in product topology. This may cause standard theory of infinite-volume states, including the variational principle, to fail, see the examples in~\cite{kulske2004relative}. 

A variety of examples have been studied ever since, where non-Gibbsian behavior was proved to occur with different mechanisms, but always in regimes of sufficiently strong coupling, where the first-layer measure differs much from independence. (Having said this, there are known examples that show that the range of temperatures where non-Gibbsian behavior in the image system occurs, may be larger than the critical temperature for the first-layer system~\cite{haggstrom2003gibbs}) Moreover, examples have been found where, the set of discontinuity points is even of full measure w.r.t.~$\mu'$ itself, which is the strongest form of singularity~\cite{kulske2004relative,EnErIaKu12,JaKu16,bergmann2020dynamical}.

Main subclasses of relevant transformations which have been studied were projections in terms of a variety of deterministic maps~\cite{van2017decimation,le2013almost,kulske2017fuzzy,MR3648046,MR3961240}, and stochastic time-evolutions~\cite{le2002short,van2002possible,kulske2007spin,EnErIaKu12,roelly2013propagation,FeHoMa14,JaKu16,kraaij2021hamilton}, in various underlying geometries of lattice models, mean-field models, Kac-models, and models in the continuum~\cite{JaKu16}. 

\subsubsection*{Informal result: Even independent fields may become non-Gibbs under projections}

In the present paper we provide a new and simple example that shows that a natural local transform of range $1$ can produce singularities, even when it is applied to an {\em independent field}. In our example we chose as the first-layer field the i.i.d.~Bernoulli lattice field $\mu_p$ on the integer lattice, with state space $\{0,1\}^{Z^d}$, and occupation probability $p\in [0,1)$. The Bernoulli lattice field in itself is studied in site-percolation, where one asks for existence of infinite clusters and refined connectedness properties~\cite{grimmett1999percolation}. It also drives more complex processes, in statistical mechanics of disordered systems~\cite{grimmett1997percolation,MR1766342,MR2252929}, and elsewhere in probability
~\cite{adler1991bootstrap,MR2283880,MR3161674,MR3156983,jahnel2020probabilistic} 
and its application. 
We then study the second-layer measure $\mu_p'$ that appears as an image under application of the concrete range-one map $T$ that is defined by removing from a realization of occupied sites the {\em occupied isolated sites}. $T$ is a projection map as it satisfies $T^2=T$, and we will call it {\em the projection to non-isolates}. Hence it keeps from a realization of occupied sites only the occupied clusters of size of at least two. This includes the infinite cluster, in case there is one, i.e., in the percolation regime of large enough $p$. We may also view $T$ as a simple smoothing transformation, as isolated 'dust' of occupied sites (or 'pixels') is forgotten under the map. 

What to expect for the second-layer measure $\mu_p'$? As the $\mu_p$-probability that a given site is isolated equals 
$p(1-p)^{2d}$,  the map $T$ seems non-invasive, in particular for probabilities $p$ close to $1$. In particular the removed sites do no percolate in this regime. So one might naively conjecture that the second-layer measure shall not be much affected and well-behaved, and in particular is representable as a Gibbsian distribution with quasilocal conditional probabilities.
  
As a main message of this paper, we prove that this is not the case: We show that $\mu'_p$ is spatially {\em non-Markovian} and {\em non-quasilocal}, when $p<1$ is large, see Theorem~\ref{thm_large_p}. We hope that this result is of some interest for the percolation community. We complement this result by proving regularity of the projected measure when $p\geq 0$ 
is small enough see, Theorem~\ref{thm_small_p}. This implies the existence of a Gibbs-non Gibbs transition driven by $p$.  

Let us finally discuss our map $T$ that projects to non-isolates (and removes the isolates) from a dual perspective.  Namely, $T$ has a natural companion map $T^*$ that is again a projection map. $T^*$ does precisely the opposite, it projects to the isolated sites (and removes the non-isolates). 

There is independent interest in the action of $T^*$ to the iid Bernoulli lattice field, for the reason that it produces the {\em thinned Bernoulli lattice field}, in which all occupied sites are separated. This thinned Bernoulli lattice field 
is relevant, as it is the lattice analogue of the well-known and much studied Mat\'ern process in the continuum~\cite{matern1960spatial,moller2010perfect,baccelli2012extremal}. The latter by definition is derived from a first-layer Poisson process in Euclidean space by removing all points in the realization that have at least one point in the Euclidean ball of radius one. 

Clearly, the second-layer measures of both maps $T,T^*$  acting jointly on the same Bernoulli lattice-field realization appear in a natural coupling. As a first guess one may conjecture from this that, either both second-layer measures are Gibbs, or both are non-Gibbs. We warn the reader that this is too naive, not only on the level of proofs, but also the statement may be false.  We leave the analysis of the companion process, the thinned Bernoulli lattice field, 
to another study. 

The paper is organized as follows. In Section~\ref{Sec_Setting} we present the setting and our main results and non-Gibbsianness and Gibbsianness for the Bernoulli field under the removal-of-isolates transformation. In Section~\ref{Sec_Setting} we present the corresponding proofs. 

\section{Setting and main results}\label{Sec_Setting}
\begin{figure}[!htpb]
	\centering
\begin{subfigure}{0.45\textwidth}
	\input{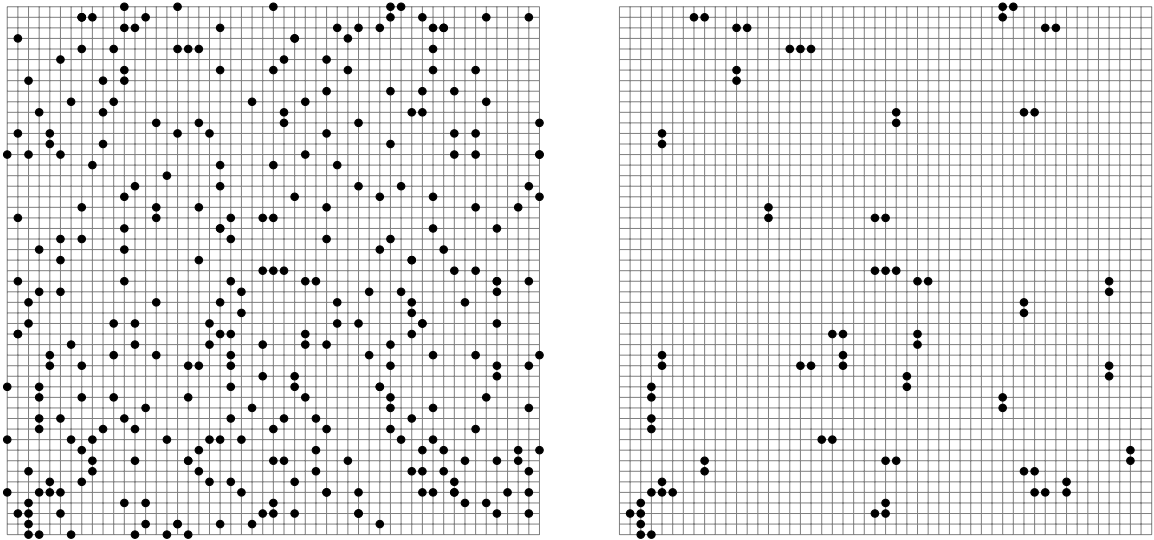}
\end{subfigure}
\begin{subfigure}{0.45\textwidth}
	\begin{tikzpicture}[scale=1.4]
\draw[step=1mm,help lines] (0,0) grid (50mm,50mm);
\filldraw[black] ( 0.10 , 0.20 ) circle (1pt);

\filldraw[black] ( 0.20 , 0.00 ) circle (1pt);
\filldraw[black] ( 0.20 , 0.10 ) circle (1pt);
\filldraw[black] ( 0.20 , 0.20 ) circle (1pt);
\filldraw[black] ( 0.20 , 0.30 ) circle (1pt);

\filldraw[black] ( 0.30 , 0.00 ) circle (1pt);
\filldraw[black] ( 0.30 , 0.40 ) circle (1pt);
\filldraw[black] ( 0.30 , 1.00 ) circle (1pt);
\filldraw[black] ( 0.30 , 1.10 ) circle (1pt);
\filldraw[black] ( 0.30 , 1.30 ) circle (1pt);
\filldraw[black] ( 0.30 , 1.40 ) circle (1pt);

\filldraw[black] ( 0.40 , 0.40 ) circle (1pt);
\filldraw[black] ( 0.40 , 0.50 ) circle (1pt);
\filldraw[black] ( 0.40 , 1.60 ) circle (1pt);
\filldraw[black] ( 0.40 , 1.70 ) circle (1pt);
\filldraw[black] ( 0.40 , 3.70 ) circle (1pt);
\filldraw[black] ( 0.40 , 3.80 ) circle (1pt);

\filldraw[black] ( 0.50 , 0.40 ) circle (1pt);

\filldraw[black] ( 0.70 , 4.90 ) circle (1pt);

\filldraw[black] ( 0.80 , 0.60 ) circle (1pt);
\filldraw[black] ( 0.80 , 0.70 ) circle (1pt);
\filldraw[black] ( 0.80 , 4.90 ) circle (1pt);

\filldraw[black] ( 1.10 , 4.30 ) circle (1pt);
\filldraw[black] ( 1.10 , 4.40 ) circle (1pt);
\filldraw[black] ( 1.10 , 4.80 ) circle (1pt);

\filldraw[black] ( 1.20 , 4.80 ) circle (1pt);

\filldraw[black] ( 1.40 , 3.00 ) circle (1pt);
\filldraw[black] ( 1.40 , 3.10 ) circle (1pt);

\filldraw[black] ( 1.60 , 4.60 ) circle (1pt);

\filldraw[black] ( 1.70 , 1.60 ) circle (1pt);
\filldraw[black] ( 1.70 , 4.60 ) circle (1pt);

\filldraw[black] ( 1.80 , 1.60 ) circle (1pt);
\filldraw[black] ( 1.80 , 4.60 ) circle (1pt);

\filldraw[black] ( 1.90 , 0.90 ) circle (1pt);

\filldraw[black] ( 2.00 , 0.90 ) circle (1pt);
\filldraw[black] ( 2.00 , 1.90 ) circle (1pt);

\filldraw[black] ( 2.10 , 1.60 ) circle (1pt);
\filldraw[black] ( 2.10 , 1.70 ) circle (1pt);
\filldraw[black] ( 2.10 , 1.90 ) circle (1pt);

\filldraw[black] ( 2.40 , 0.20 ) circle (1pt);
\filldraw[black] ( 2.40 , 2.50 ) circle (1pt);
\filldraw[black] ( 2.40 , 3.00 ) circle (1pt);

\filldraw[black] ( 2.50 , 0.20 ) circle (1pt);
\filldraw[black] ( 2.50 , 0.30 ) circle (1pt);
\filldraw[black] ( 2.50 , 0.70 ) circle (1pt);
\filldraw[black] ( 2.50 , 2.50 ) circle (1pt);
\filldraw[black] ( 2.50 , 3.00 ) circle (1pt);

\filldraw[black] ( 2.60 , 0.70 ) circle (1pt);
\filldraw[black] ( 2.60 , 2.50 ) circle (1pt);
\filldraw[black] ( 2.60 , 3.90 ) circle (1pt);
\filldraw[black] ( 2.60 , 4.00 ) circle (1pt);

\filldraw[black] ( 2.70 , 1.40 ) circle (1pt);
\filldraw[black] ( 2.70 , 1.50 ) circle (1pt);

\filldraw[black] ( 2.80 , 1.80 ) circle (1pt);
\filldraw[black] ( 2.80 , 1.90 ) circle (1pt);
\filldraw[black] ( 2.80 , 2.40 ) circle (1pt);

\filldraw[black] ( 2.90 , 2.40 ) circle (1pt);

\filldraw[black] ( 3.60 , 1.20 ) circle (1pt);
\filldraw[black] ( 3.60 , 1.30 ) circle (1pt);
\filldraw[black] ( 3.60 , 4.90 ) circle (1pt);
\filldraw[black] ( 3.60 , 5.00 ) circle (1pt);

\filldraw[black] ( 3.70 , 5.00 ) circle (1pt);

\filldraw[black] ( 3.80 , 0.60 ) circle (1pt);
\filldraw[black] ( 3.80 , 2.10 ) circle (1pt);
\filldraw[black] ( 3.80 , 2.20 ) circle (1pt);
\filldraw[black] ( 3.80 , 4.00 ) circle (1pt);

\filldraw[black] ( 3.90 , 0.40 ) circle (1pt);
\filldraw[black] ( 3.90 , 0.60 ) circle (1pt);
\filldraw[black] ( 3.90 , 4.00 ) circle (1pt);

\filldraw[black] ( 4.00 , 0.40 ) circle (1pt);
\filldraw[black] ( 4.00 , 4.80 ) circle (1pt);

\filldraw[black] ( 4.10 , 4.80 ) circle (1pt);

\filldraw[black] ( 4.20 , 0.40 ) circle (1pt);
\filldraw[black] ( 4.20 , 0.50 ) circle (1pt);

\filldraw[black] ( 4.60 , 1.50 ) circle (1pt);
\filldraw[black] ( 4.60 , 1.60 ) circle (1pt);
\filldraw[black] ( 4.60 , 2.30 ) circle (1pt);
\filldraw[black] ( 4.60 , 2.40 ) circle (1pt);

\filldraw[black] ( 4.80 , 0.70 ) circle (1pt);
\filldraw[black] ( 4.80 , 0.80 ) circle (1pt);
\end{tikzpicture}
\end{subfigure}
	\caption{Realization of a Bernoulli field (left) and its image under the transformation $T$ (right).}
	\label{Pix_Trans}
\end{figure}

To define our process we start from the $\Omega=\{0,1\}^{\Z^d}$-valued i.i.d.~Bernoulli field $\mu_p$ with parameter $p\in [0,1]$. We consider realizations of the Bernoulli field under the application of the transformation 
$T:\Omega\rightarrow \Omega$ 
given by 
$$
(T\omega)_x:=\o'_x=\omega_x\Bigl(1-\prod_{y\in \partial x}(1-\omega_y)\Bigr),
$$
where $\partial x$ denotes the set of nearest neighbors of $x\in\Z^d$ in $\Z^d$, equipped with the usual neighborhood structure. In words, $T$ is the projection to the non-isolates, see Figure~\ref{Pix_Trans} for an illustration. The image measure under the transformation $$\mu'_p:=\mu_p\circ T^{-1}$$ is supported on the subset $\O':=T(\O)$ of sites that obey the non-isolation constraint. 

Intuitively the application of $T$ should not change the measure very much at large $p$ close to $1$, where a typical configuration consists of a large percolating cluster and very few isolated sites. In this regime one may view $T$ as a cleansing transformation that wipes away the smallest dust of isolated sites, and keeps the apparent main parts. From the definition of $T$ as a local map it is also obvious that for variables at sites of graph distance greater equal than $3$ are independent, under $\mu'_p$, so that one could expect that $\mu'_p$ is a nicely behaved measure. 

Recall that a {\em specification} $\g=(\g_\L)_{\L\Subset\Z^d}$ is a consistent and proper family of conditional probabilities, i.e., for all $\L\subset\D\Subset\Z^d$, $\o_\L\in \O_\L:=\{0,1\}^\L$ and $\hat\o\in \O$, we have that $\int_{\O}\g_\D(\d \tilde\o|\hat\o)\g_\L(\o_\L|\tilde\o)=\g_\D(\o_\L|\hat\o)$, and for all $\o_{\L^{\rm c}}\in \O_{\L^{\rm c}}$ we have $\g_\L(\o_{\L^{\rm c}}|\hat\o)=\one\{\o_{\L^{\rm c}}=\hat\o_{\L^{\rm c}}\}$. A specification is called {\em quasilocal}, if for all $\L\Subset\Z^d$ and $\hat\o_\L \in\O_\L$, the mapping $\o\mapsto\g_\L(\hat\o_\L|\o)$ is continuous with respect to the product topology on $\O$. We say that $\g$ is a specification for some random field $\mu$ on $\O$, if it satisfies the DLR equations, i.e., for all $\L\Subset\Z^d$ and $\o_\L\in \O_\L$, we have that $\int_{\O}\mu(\d \tilde\o)\g_\L(\o_\L|\tilde\o)=\mu(\o_\L)$.

Our present result shows that this is not the case in the whole parameter regime, but the following is true: For large $p$, $\mu'_p$ is not spatially Markovian, it is not even a Gibbs measure in the sense of existence of a version of its finite-volume conditional probabilities which is continuous with respect to the product topology on $\O'$. More precisely we have the following theorem. 
\begin{thm}(Non-Gibbsianness for large $p$)\label{thm_large_p} Consider the image measure $\mu'_p$ of the Bernoulli field on $\Z^d$ under the map to the non-isolates in lattice dimensions $d\geq 2$. Then, there is $p_c(d)<1$ such that for $p\in (p_c(d),1)$, there is no quasilocal specification 
$\g'$ for $\mu'_p$. 
\end{thm}

The result shows that Gibbsian descriptions of thinning processes of various types derived from Bernoulli or Poissonian fields are by no means obvious, and that more research on such processes in discrete and continuous setups is necessary. 

\begin{figure}[!htpb]
\centering
\begin{tikzpicture}[scale=1]
\draw (-0.1,0) -- (10,0);
\draw (-0.1,1) -- (0,1);
\draw (-0.1,2) -- (0,2);
\draw (-0.1,3) -- (0,3);
\draw (0,0) -- (0,3);
\draw [fill=green,opacity=0.3](0,0) rectangle (10,1);
\draw [fill=red,opacity=0.3](0,2) rectangle (10,3);
\draw [fill=green,opacity=0.3](0,2.98) rectangle (10,3);
\node at (-0.3,3) {$1$};
\node at (-0.6,2) {$p_c(d)$};
\node at (-0.6,1) {$\frac{1}{2d}$};
\node at (-0.3,0) {$0$};
\node at (5,0.5) {Gibbsianness};
\node at (5,2.5) {non-Gibbsianness};
\end{tikzpicture}
\caption{Illustration of Gibbs-non-Gibbs transitions in $p$ for the thinned Bernoulli field under non-isolation constraint.}
\label{fig_1}
\end{figure}
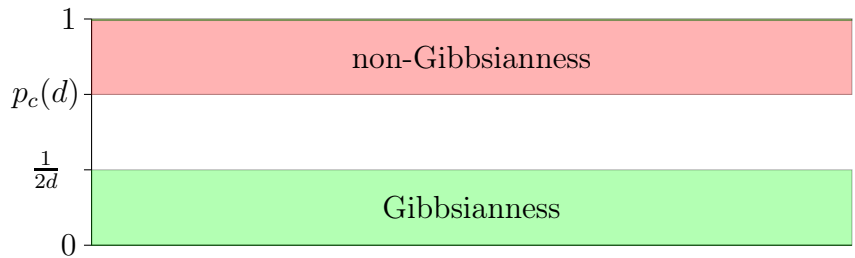

In order to complement the above non-Gibbsianness result, let us also present the following theorem on the existence of a quasilocal specification for $\mu'_p$ for small values of $p$. 
\begin{thm}(Gibbsianness for small $p$)\label{thm_small_p} Consider the image measure $\mu'_p$ of the Bernoulli field on $\Z^d$ under the map to the non-isolates in lattice dimensions $d\geq 1$. Then,
for $p< 1/(2d)$, there exists a quasilocal specification $\gamma'$ for $\mu'_p$. 
\end{thm}
As we can see from the proofs, for $p< 1/(2d)$, in fact the continuity of $\g'$ is even exponentially fast. 
In summary, the statements of Theorems~\ref{thm_large_p} and \ref{thm_large_p} indicate a phase-diagram of Gibbsianness of thinned Bernoulli fields under the local non-isolation constraint as exhibited in Figure~\ref{fig_1}.

In the following section, we present the proofs. 

\section{Proofs}\label{Sec_Proofs}
The key idea of the proof is to first re-express conditional probabilities of $\mu'_p$ in finite volumes in terms of a first-layer constraint model in which occupied sites have to be isolated. Indeed, for large $p$, there are two distinct groundstates given by the (shifted) checkerboard configurations. We can leverage a Peierls' argument in order to show that the first-layer constraint model exhibits a phase transition of translational symmetry breaking. We note that the argument works even though there is {\em no} spin-flip symmetry in the system. The translational symmetry breaking gives rise to a point of discontinuity for which we subsequently show that it is present for any system of finite-volume conditional probabilities for $\mu'_p$, i.e., it is an essential discontinuity. The converse case of small $p$ can be handled by arguments using Dobrushin uniqueness techniques. 

\subsection{Proof of Theorem~\ref{thm_large_p}: Non-Gibbianness}
\label{Sec_Proofs_High_p_s}
The main ingredient for the proof is to exhibit one non-removable bad configuration for conditional probabilities. For this, we will use the so-called two-layer view, in which one needs to understand 
the Bernoulli field conditional on a fixed image configuration. We choose as the image configuration the all empty configuration, for which the first-layer measure becomes 
the Bernoulli field $\mu_p$ {\em conditional on isolation}. 

We proceed as follows. In Section~\ref{Sec_Pei} we exhibit a phase transition for the latter model at large $p$, in which translation symmetry is broken, which can be selected via suitable shapes of loophole-volumes. The technique is based on a (slightly non-standard) Peierls argument.  
In Section~\ref{Sec_Non} we then show how this implies non-Gibbsianness of the image measure. This is based on the proof that jumps of conditional probabilities occur for certain suitably chosen local patterns, which allow to make a transparent connection to the first-layer model in suitable connected boxes, where the Peierls argument from Section~\ref{Sec_Pei} was made to work. 
 
\subsubsection{Translational-symmetry breaking via a Peierls argument for the conditional first-layer model}\label{Sec_Pei}
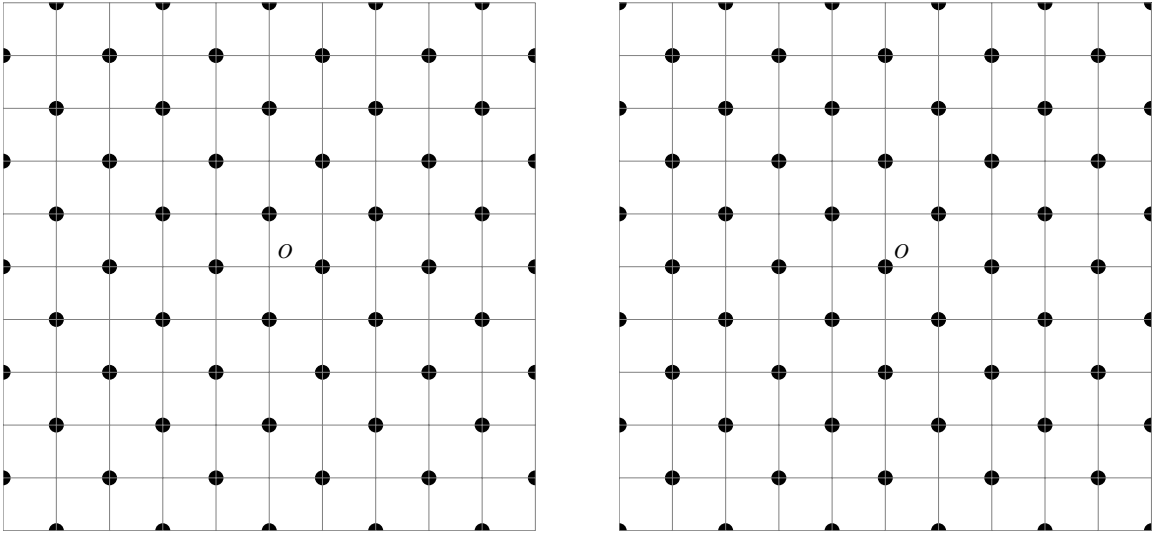
\begin{figure}[!htpb]
	\centering
\begin{subfigure}{0.45\textwidth}
	\begin{tikzpicture}[scale=1.4]
\begin{scope}
\clip (-25mm,-25mm) rectangle (25mm,25mm);
\foreach \i in {-3,...,2}
{
\foreach \j in {-2,...,2}
{
\fill (\i+0.5,\j)  circle[radius=2pt];
}}
\foreach \i in {-2,...,2}
{
\foreach \j in {-3,...,2}
{
\fill (\i,\j+0.5)  circle[radius=2pt];
}}
\node at (0.15,0.15) {$o$};
\draw[step=5mm,help lines] (-25mm,-25mm) grid (25mm,25mm);
\end{scope}
\end{tikzpicture}
\end{subfigure}
\begin{subfigure}{0.45\textwidth}
	\begin{tikzpicture}[scale=1.4]
\begin{scope}
\clip (-25mm,-25mm) rectangle (25mm,25mm);
\foreach \i in {-2,...,2}
{
\foreach \j in {-2,...,2}
{
\fill (\i,\j)  circle[radius=2pt];
}}
\foreach \i in {-3,...,2}
{
\foreach \j in {-3,...,2}
{
\fill (\i+0.5,\j+0.5)  circle[radius=2pt];
}}
\node at (0.15,0.15) {$o$};
\draw[step=5mm,help lines] (-25mm,-25mm) grid (25mm,25mm);
\end{scope}
\end{tikzpicture}
\end{subfigure}
	\caption{Illustration of the type-0 (left) and type-1 (right) groundstates, i.e., checkerboard configurations. Dots indicate occupied sites. The origin is indicated as $o$. }
	\label{Pix_Types}
\end{figure}

For the purpose of showing that the empty configuration is bad for any specification, we will analyze the following particular finite-volume first-layer measures, and we will restrict to particular volumes $\L$. Namely, let us consider finite volumes $\Lambda$ that have a {\em shape of type $0$}, and put fully occupied boundary conditions all $1$ outside of $\Lambda$. By this we mean that $\Lambda$ has a shape which allows to put the checkerboard groundstate of zeros and ones inside $\Lambda$ for which the origin obtains the value $0$ such that one obtains a configuration compatible with the boundary condition, see Figure~\ref{Pix_Types} for an illustration. We define
\begin{equation*}
\begin{split}
\nu_{\L}( \omega_{\L} ):= 
\frac{\mu_{p,\Lambda}( \omega_{\L}  1_{T(\omega_{\L}1_{\Lambda^c})|_{\L}=0_{\Lambda}} )}{
\mu_{p,\Lambda}(T(\sigma_{\L}1_{\Lambda^c})|_{\L}=0_{\Lambda} )
},
\end{split}
\end{equation*}
where $\mu_{p,\Lambda}$ is the Bernoulli product measure in $\Lambda$. Hence, by definition $\nu_{\Lambda}$ is the Bernoulli measure conditioned on isolation of ones inside $\Lambda$, where the isolation constraint remembers also the fully occupied boundary condition.

\begin{figure}[!htpb]
	\centering
	\begin{tikzpicture}[scale=1.4]
\draw[step=1mm,help lines] (-25mm,-25mm) grid (25mm,25mm);

\foreach \i in {-9,...,10}
{
\filldraw[black] (\i*0.2-0.1,1.8) circle (1pt);
}

\foreach \j in {0,...,6}
{\foreach \i in {-25,...,25}
{
\filldraw[black] (\i*0.1,1.9+\j*0.1) circle (1pt);
}
}

\foreach \j in {0,...,6}
{\foreach \i in {-25,...,25}
{
\filldraw[black] (\i*0.1,-1.9-\j*0.1) circle (1pt);
}
}

\foreach \j in {0,...,6}
{\foreach \i in {-25,...,25}
{
\filldraw[black] (1.9+\j*0.1,\i*0.1) circle (1pt);
}
}

\foreach \j in {0,...,6}
{\foreach \i in {-25,...,25}
{
\filldraw[black] (-1.9-\j*0.1,\i*0.1) circle (1pt);
}
}

\foreach \i in {-9,...,10}
{
\filldraw[black] (\i*0.2-0.1,-1.8) circle (1pt);
}
\foreach \i in {-9,...,9}
{
\filldraw[black] (\i*0.2,-1.9) circle (1pt);
}

\foreach \i in {-9,...,9}
{
\filldraw[black] (1.8,\i*0.2-0.1) circle (1pt);
}
\foreach \i in {-9,...,9}
{
\filldraw[black] (1.9,\i*0.2) circle (1pt);
}
\foreach \i in {-9,...,9}
{
\filldraw[black] (-1.8,\i*0.2-0.1) circle (1pt);
}
\foreach \i in {-9,...,9}
{
\filldraw[black] (-1.9,\i*0.2) circle (1pt);
}

\foreach \i in {-8,...,8}
{
\draw[blue] (\i*0.2-0.05,1.85) -- (\i*0.2+0.05,1.85);
}

\foreach \i in {-9,...,8}
{
\draw[blue] (\i*0.2+0.05,1.75) -- (\i*0.2+0.15,1.75);
}

\foreach \i in {-18,...,19}
{
\draw[blue] (\i*0.1-0.05,1.85) -- (\i*0.1-0.05,1.75);
}

\foreach \i in {-8,...,8}
{
\draw[blue] (\i*0.2-0.05,-1.85) -- (\i*0.2+0.05,-1.85);
}

\foreach \i in {-9,...,8}
{
\draw[blue] (\i*0.2+0.05,-1.75) -- (\i*0.2+0.15,-1.75);
}

\foreach \i in {-18,...,19}
{
\draw[blue] (\i*0.1-0.05,-1.85) -- (\i*0.1-0.05,-1.75);
}

\foreach \i in {-8,...,8}
{
\draw[blue] (1.85,\i*0.2-0.05) -- (1.85,\i*0.2+0.05);
}

\foreach \i in {-9,...,8}
{
\draw[blue] (1.75,\i*0.2+0.05) -- (1.75,\i*0.2+0.15);
}

\foreach \i in {-18,...,19}
{
\draw[blue] (1.85,\i*0.1-0.05) -- (1.75,\i*0.1-0.05);
}

\foreach \i in {-8,...,8}
{
\draw[blue] (-1.85,\i*0.2-0.05) -- (-1.85,\i*0.2+0.05);
}

\foreach \i in {-9,...,8}
{
\draw[blue] (-1.75,\i*0.2+0.05) -- (-1.75,\i*0.2+0.15);
}

\foreach \i in {-18,...,19}
{
\draw[blue] (-1.85,\i*0.1-0.05) -- (-1.75,\i*0.1-0.05);
}

\end{tikzpicture}
	\caption{Illustration of a type-0 volume with loophole boundary. 
	}
	\label{Pix_Loop}
\end{figure}
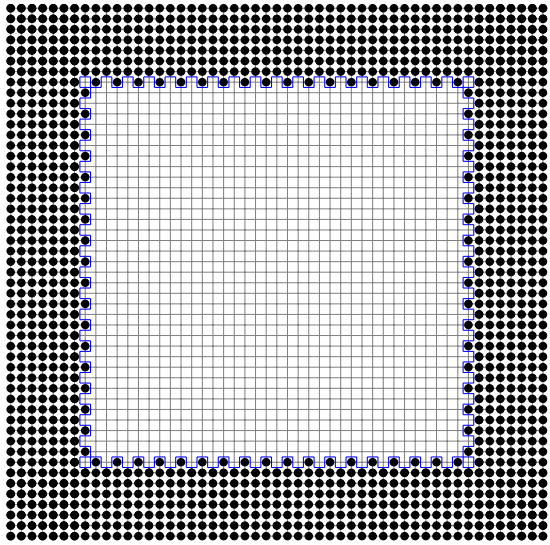

A similar definition is made for volumes of type $1$. For us, large boxes $B_L$ centered around the origin with sidelength $2L$ with a loophole boundary, will be useful, see the illustration in Figure~\ref{Pix_Loop}.  
For such type-$0$ boxes $B_L$ we will show in this section that 
\begin{equation}\label{eq_0}
\begin{split}
&\sup_{L}\nu_{B_L}(\omega_0=1)\leq \epsilon(p),\qquad \text{ with }\qquad\lim_{p\uparrow 1}\epsilon(p)=0.
\end{split}
\end{equation}
This means that, with large probability, the origin copies the information from the boundary condition. Similarly, we will prove that the spin at the origin for the box shifted by a lattice unit vector $e$ satisfies 
\begin{equation}\label{eq_1}
\begin{split}
&\sup_{L}\nu_{B_L+e}(\omega_0=0)\leq \epsilon(p),\qquad \text{ with }\qquad\lim_{p\uparrow 1}\epsilon(p)=0.
\end{split}
\end{equation}
This is an essential step as it proves that the shape of the volume $B_L$ induces a phase transition for the first-layer constrained model, and there is breaking of translational symmetry. To complete the proof of essential badness of the empty configuration on the second layer, we will however need to go one step further, and connect 
to the measure on the second layer. This will be done in Lemma~\ref{lem_NonGibbs} below. 

Note that configurations of the model are energetically equivalent under a lattice shift.  They are not equivalent under the site-wise {\em spin-flip} that exchanges zeros and ones, much unlike the Ising antiferromagnet in zero external field. Therefore, the Peierls argument we are about to give has to be different from the one for the Ising ferromagnet or antiferromagnet. 
Namely, the Peierls argument we will present involves suitable lattice shifts of parts of configurations, while the standard more straightforward Peierls argument for the Ising model involves spin-flips.   

Consider the nearest-neighbor graph with vertex set $\Z^d$. Consider, for a spin configuration $\omega$, the set of sites 
$$\Gamma(\omega):=\{x \in \Z^d\colon \text{ there exists } y\in \partial x \text{ such that }\omega_x=\omega_y=0 \}.$$
Note that there is a one-to-one correspondence between configurations $\omega$ that satisfy the neighborhood constraint, and sets $\Gamma(\omega)$. Note that outside of $\Gamma(\omega)$, the configuration $\omega$ looks like one of the two groundstates formed by the two possible checkerboard configurations of zeros and ones. Indeed, each site $x\not\in \Gamma(\omega)$ has the property that either $\omega_x=0$ and all the neighbors are $1$ (by definition of a contour), 
or $\omega_x=1$ and all the neighbors are $0$ (as the model contains the isolation-constraint).  

Further note that not all possible subsets of $\Z^d$ can occur as $\Gamma$, because of the isolation constraint of ones.  

\begin{figure}[!htpb]
	\centering
	\begin{tikzpicture}[scale=1.4]
\begin{scope}
\clip (-15mm,-15mm) rectangle (15mm,15mm);
\draw[step=5mm,help lines] (-25mm,-25mm) grid (25mm,25mm);
\foreach \i in {-2,...,2}
{
\foreach \j in {-2,...,2}
{
\fill (\i,\j)  circle[radius=2pt];
}}
\foreach \i in {-3,...,2}
{
\foreach \j in {-3,...,2}
{
\fill (\i+0.5,\j+0.5)  circle[radius=2pt];
}}
\node at (0.15,0.15) {$o$};
\fill[white] (0,0)  circle[radius=3pt];
\draw[step=5mm,help lines] (-25mm,-25mm) grid (25mm,25mm);
\draw[blue] (-0.25,0+0.75) -- (0.25,0+0.75);
\draw[blue] (-0.25,0-0.75) -- (0.25,0-0.75);
\draw[blue] (-0.75,0+0.25) -- (-0.25,0+0.25);
\draw[blue] (0.25,0+0.25) -- (0.75,0+0.25);
\draw[blue] (-0.75,0-0.25) -- (-0.25,0-0.25);
\draw[blue] (0.25,0-0.25) -- (0.75,0-0.25);

\draw[blue] (0+0.75,-0.25) -- (0+0.75,0.25);
\draw[blue] (0-0.75,-0.25) -- (0-0.75,0.25);
\draw[blue] (0+0.25,-0.75) -- (0+0.25,-0.25);
\draw[blue] (0+0.25,0.25) -- (0+0.25,0.75);
\draw[blue] (0-0.25,-0.75) -- (0-0.25,-0.25);
\draw[blue] (0-0.25,0.25) -- (0-0.25,0.75);

\end{scope}
\end{tikzpicture}
	\caption{Illustration of a minimal contour in blue containing the origin. }
	\label{Pix_MinCont}
\end{figure}
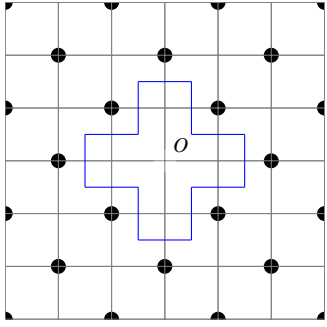

The connected (in the sense of graph-distance) components $\gamma$ of these sets $\Gamma$ are called {\em contours}. 
To visualize this, consider a star-shaped contour that is built from flipping the one site from one to zero starting from 
a checkerboard configuration, see Figure~\ref{Pix_MinCont}. This yields the minimal contour which has $2d+1$ sites. In two dimensions, e.g., it is possible that different contours can be reached from each other via the diagonal. Note that each $\gamma$ that is a contour of a configuration, must be surrounded by ones in nearest neighbor sense in the configuration. These ones must be surrounded by nearest neighbors which carry zeros, by the isolation constraint of ones. Hence the contour specifies the configuration up to sites with graph distance two. 
 
The complement of a finite contour $\gamma$ has one infinite component, and finitely many finite components (the internal ones). Each of these components are labelled by one of the two labels $1$ (or $0$ respectively) determined whether, given $\gamma$, the component admits a configuration obtained by substituting the infinite-volume checkerboard configurations in which the origin in $\Z^d$ obtains a $1$ (or a $0$ respectively), see Figure~\ref{Pix_Types}. 
A configuration is uniquely determined by its set of compatible contours. Contours $\gamma$ themselves are labelled by $1$ (or $0$) according to the two possible labels of their outer connected components. Contours are compatible when they arise from an allowed configuration, which is the case when the types of checkerboards on shared connected components of the complements match. 

Now, suppose that $\gamma$ is a contour in a finite volume $\Lambda$. We decompose the volume in terms of the contour and the connected components of its complement 
$$\Lambda=\gamma\cup V_0\cup_{i=1,\dots,k}V_i.$$ 
Here $V_0$ is the outer connected component of the complement of $\gamma$ intersected with $\Lambda$. (Note that the intersection with $\Lambda$ may produce several connected components, but this poses no difficulty.)  The sets $(V_i)_{i\geq 1}$ are the interior connected components of the complement of $\gamma$. We also write  
$$\gi=\gamma\cup_{i=1,\dots,k}V_i$$
for the sites that are contained in the support of $\gamma$ or surrounded by $\gamma$.  

\medskip
We now prove the statements~\eqref{eq_0} and~\eqref{eq_1} via a Peierls estimate. It suffices 
to treat the first case~\eqref{eq_0}, as the case~\eqref{eq_1} is similar.  We start with a union bound 
over contours surrounding the origin
\begin{equation*}
\begin{split}
\nu_{B_L}(\omega_0=1)&\leq\nu_{B_L}(\omega\colon \exists \gamma \text{  such that } 
\Gamma(\omega)\ni \gamma \text{ and } \gi \ni 0)\leq \sum_{\gamma\colon   \gi \ni 0 } 
\nu_{B_L}(\omega\colon\Gamma(\omega)\ni \gamma ).
\end{split}
\end{equation*}
With a slight abuse of notation, we here write $\gamma \in \Gamma$ to indicate that $\gamma$ is a contour in $\Gamma$. 

The main point for the Peierls estimate in our non-flip invariant situation, which is formulated in the following lemma, 
will be the construction of compatible configurations after removal of contours. Let $|\gamma|$ denote the number of vertices in $\gamma$. 
\begin{lem}\label{lem_Pei}
There exists a Peierls constant $\tau=\tau(p)$ with $\lim_{p\uparrow 1}\tau(p)=\infty$, such that for all $\gamma$ we have that
\begin{equation*}
\begin{split}
&\nu_{B_L}(\omega\colon \Gamma(\omega)\ni \gamma)\leq e^{-\tau |\gamma|}.
\end{split}
\end{equation*}
\end{lem}
\begin{proof}Define the activity of a contour $\gamma\subset B_L$ 
to be the natural weight of the zeros prescribed by it in the Bernoulli measure, i.e., 
$$
\rho(\gamma):=(1-p)^{|\gamma|}.
$$
For any configuration $\omega_U$ in a finite volume $U$ 
we write for its weight in the Bernoulli field 
$$
R(\omega_U):=\prod_{x\in U }p^{\omega_x}(1-p)^{1-\omega_x}.
$$
In particular $\rho(\gamma)=R(0_\gamma)$, where we use the short-hand notation $0_B$ to indicate the all-zero configuration in the volume $B$. Then, we may write 
\begin{equation*}
\begin{split}
&\nu_{B_L}(\omega: \Gamma(\omega)\ni \gamma )=\frac{\rho(\gamma)Z_{V_0}\prod_{i=1}^k Z_{V_i}}{Z_{B_L}},
\end{split}
\end{equation*}
where $Z_{V_i}$, for $i=0,1,\dots,k$ denotes the partition functions over all configurations in the volumes $V_i$ which are compatible with $\gamma$, with the isolation constraint on the ones, with weights provided by the Bernoulli measure for all sites in $V_i$.

\paragraph{Case 1.} Suppose that $\gamma$ only contains interior connected components of type $0$. This means that there is no typechange when going from the outside to the inside. Then we may remove $\gamma$, i.e., continue the checkerboard configuration outside of $\gamma$ to where $\gamma$ used to be. 
This means that we will assign to each $\omega$ for which $\Gamma(\omega)\ni \gamma$ the reference 
configuration $(\omega_{\Lambda \backslash\gamma}\omega^0_{\gamma})$ which appears in the partition function $Z_{\Lambda}$ to lower bound the latter. Here $\omega^0_{\gamma}$ denotes the type-0 checkerboard configuration on $\gamma$. 
It is important to note that this removal keeps all other contributions from exterior and interior components compatible. 
Hence, we immediately arrive at a lower bound 
\begin{equation*}
\begin{split}
&Z_{B_L}\geq R(\omega^0_\gamma)Z_{V_0}\prod_{i=1}^k Z_{V_i}.
\end{split}
\end{equation*}
We may write $\rho(\gamma)=R(\omega^0_\gamma)((1-p)/p)^{N^\text{repl}}$, where $N^\text{repl}$ denotes the number of replacements of a zero by a one on the support of the contour. By definition of a contour each site in $\gamma$ has a neighbor which is zero. Therefore it will replaced itself by a one, or a neighbor of it will be replaced by a one. Hence $N^\text{repl}\geq|\gamma|/(2d +1)$. Thus, we have the desired estimate 
\begin{equation*}
\begin{split}
&\nu_{B_L}(\omega: \Gamma(\omega)\ni \gamma )\leq \big((1-p)/p\big)^{|\gamma|/(2d +1)}.
\end{split}
\end{equation*}

\paragraph{Case 2.} Suppose now that $\gamma$ additionally also contains interior volumes of type $1$ (the bad type that does not agree with the boundary outside $\L$), which we will denote by $W_j$, $j=1,\dots,l$. Writing $V_i$ for the interior components of type-0, we have
\begin{equation*}
\begin{split}
&\nu_{B_L}(\omega: \Gamma(\omega)\ni \gamma )=\frac{\rho(\gamma)Z_{V_0}\prod_{i=l+1}^k Z_{V_i}\prod_{j=1}^l Z_{W_j}}{Z_{B_L}},
\end{split}
\end{equation*}
where all partition functions are sums over compatible configurations in the respective connected components, such that the total configuration contains the contour $\gamma$.

The difficulty of this case is that the removal of $\gamma$ does not immediately create compatible configurations. However, it does so after the shift of each of the internal volumes of wrong checkerboard subtypes $W_j$ in one of the $2d$ possible (positive or negative) lattice directions $e$. Let us explain the details now. 
Our comparison configuration will now be equal to the type-$0$ checkerboard on the following set $\gamma_e$ which describes the appropriate modification of $\gamma$, obtained by shifts of the internal components, 
\begin{equation}\label{Con_Mov}
\begin{split}
\gamma_e:=\Bigl(\gamma \backslash\bigcup_{j=1}^l(W_j+e) \Bigr)\cup 
\bigcup_{j=1}^lW_j \backslash (W_j+e).
\end{split}
\end{equation}
Note that there is the volume preservation $|\gamma_e|=|\gamma|$. For each $\omega$ for which $\Gamma(\omega)\ni \gamma$, the reference configuration will then be 
\begin{equation}\label{Ref_Con}
\begin{split}
(\omega^0_{\gamma_e},\omega_{\cup_{i=l+1}^k V_i},(\theta_e \omega)_{\cup_{j=1}^l (W_j+e)}),
\end{split}
\end{equation}
where $\theta_e$ represents the shift of the configuration by $e$. We see from the definition that $\omega$ will not be modified on the external component and the internal components of good type. It will however be shifted by $e$ on the internal components of bad type, and it will be a checkerboard of good type on the modified contour $\gamma_e$. Note that this configuration really occurs in $Z_{\Lambda}$ as it satisfies the isolation constraint on the ones. Hence it can be used to lower bound the partition function which gives us 
\begin{equation*}
\begin{split}
&Z_\Lambda\geq R(\omega^0_{\gamma_e})
Z_{V_0}\prod_{i=l+1}^k Z_{V_i}\prod_{j=1}^l Z_{W_j},
\end{split}
\end{equation*}
where we have used the shift invariance for the internal partition functions $Z_{W_j}$ (those with the bad types).  

Now, the proof is finished once we show that there is a dimension-dependent constant $c_d>0$ such that 
\begin{equation}\label{cd}
\begin{split}
&\rho(\gamma)\leq R(\omega^0_{\gamma_e})\big((1-p)/p\big)^{c_d |\gamma|}.
\end{split}
\end{equation}
We denote by $S^0$ the occupied sites of $\omega^0$ (the good checkerboard configuration).  

The idea is to use the fact that any connected (in graph distance) subset of $\Z^d$ hits at least a positive fraction of $S^0$ to conclude the inequality 
\begin{equation}\label{cd2}
\begin{split}
&|\gamma_e \cap S^0|\geq c_d |\gamma_e|=c_d |\gamma|.
\end{split}
\end{equation}
Then, the Inequality~\eqref{cd} would follow immediately from that.  

However, there is the small problem with that argument since, while $\gamma$ by definition is connected (w.r.t.~graph distance), the modified set $\gamma_e$ may have obtained isolated sites, see Figure~\ref{Pix_Iso}. 
\begin{figure}[!htpb]
	\centering
\begin{subfigure}{0.325\textwidth}
\begin{tikzpicture}[scale=3]
\begin{scope}
\clip (-9mm,-9mm) rectangle (9mm,9 mm);
\draw[step=1mm,help lines] (-25mm,-25mm) grid (25mm,25mm);

\foreach \i in {-4,...,5}
{
\filldraw[black] (\i*0.2-0.1,0.9) circle (1pt);
}

\foreach \i in {-4,...,5}
{
\filldraw[black] (\i*0.2,0.8) circle (1pt);
}

\foreach \i in {-4,...,5}
{
\filldraw[black] (\i*0.2-0.1,0.7) circle (1pt);
}

\foreach \i in {-4,...,5}
{
\filldraw[black] (\i*0.2,0.6) circle (1pt);
}

\foreach \i in {-4,...,5}
{
\filldraw[black] (\i*0.2-0.1,0.5) circle (1pt);
}

\foreach \i in {-1,...,1}
{
\filldraw[green] (\i*0.2,0.5) circle (1pt);
}
\foreach \i in {-1,...,1}
{
\filldraw[green] (\i*0.2,-0.5) circle (1pt);
}

\foreach \i in {-4,...,-2}
{
\filldraw[black] (\i*0.2,0.4) circle (1pt);
}

\foreach \i in {-1,...,2}
{
\filldraw[green] (\i*0.2,0.4) circle (1pt);
}
\foreach \i in {-1,...,2}
{
\filldraw[green] (\i*0.2,-0.4) circle (1pt);
}
\foreach \i in {-2,...,3}
{
\filldraw[green] (\i*0.2-0.1,0.4) circle (1pt);
}
\foreach \i in {-2,...,3}
{
\filldraw[green] (\i*0.2-0.1,-0.4) circle (1pt);
}

\foreach \i in {2,...,5}
{
\filldraw[black] (\i*0.2,0.4) circle (1pt);
}

\foreach \i in {-4,...,-3}
{
\filldraw[black] (\i*0.2-0.1,0.3) circle (1pt);
}
\foreach \i in {4,...,5}
{
\filldraw[black] (\i*0.2-0.1,0.3) circle (1pt);
}
\foreach \i in {-1,...,1}
{
\filldraw[black] (\i*0.2,0.3) circle (1pt);
}

\foreach \i in {-6,...,-3}
{
\filldraw[green] (\i*0.1,0.3) circle (1pt);
}
\foreach \i in {3,...,6}
{
\filldraw[green] (\i*0.1,0.3) circle (1pt);
}
\foreach \i in {-6,...,-3}
{
\filldraw[green] (\i*0.1,-0.3) circle (1pt);
}
\foreach \i in {3,...,6}
{
\filldraw[green] (\i*0.1,-0.3) circle (1pt);
}
\foreach \i in {-1,...,1}
{
\filldraw[green] (\i*0.2-0.1,0.3) circle (1pt);
}
\foreach \i in {-1,...,1}
{
\filldraw[green] (\i*0.2-0.1,-0.3) circle (1pt);
}

\foreach \i in {-4,...,-3}
{
\filldraw[black] (\i*0.2,0.2) circle (1pt);
}
\foreach \i in {3,...,5}
{
\filldraw[black] (\i*0.2,0.2) circle (1pt);
}
\foreach \i in {-1,...,2}
{
\filldraw[black] (\i*0.2-0.1,0.2) circle (1pt);
}

\foreach \i in {-5,...,-4}
{
\filldraw[green] (\i*0.1,0.2) circle (1pt);
}
\foreach \i in {4,...,5}
{
\filldraw[green] (\i*0.1,0.2) circle (1pt);
}
\foreach \i in {-5,...,-4}
{
\filldraw[green] (\i*0.1,-0.2) circle (1pt);
}
\foreach \i in {4,...,5}
{
\filldraw[green] (\i*0.1,-0.2) circle (1pt);
}

\foreach \i in {-4,...,-2}
{
\filldraw[black] (\i*0.2-0.1,0.1) circle (1pt);
}
\foreach \i in {3,...,5}
{
\filldraw[black] (\i*0.2-0.1,0.1) circle (1pt);
}
\foreach \i in {-1,...,1}
{
\filldraw[black] (\i*0.2,0.1) circle (1pt);
}

\foreach \i in {-4,...,-3}
{
\filldraw[green] (\i*0.1,0.1) circle (1pt);
}
\foreach \i in {3,...,4}
{
\filldraw[green] (\i*0.1,0.1) circle (1pt);
}
\foreach \i in {-4,...,-3}
{
\filldraw[green] (\i*0.1,-0.1) circle (1pt);
}
\foreach \i in {3,...,4}
{
\filldraw[green] (\i*0.1,-0.1) circle (1pt);
}

\foreach \i in {-4,...,-3}
{
\filldraw[black] (\i*0.2,0) circle (1pt);
}
\foreach \i in {3,...,5}
{
\filldraw[black] (\i*0.2,0) circle (1pt);
}
\foreach \i in {-1,...,2}
{
\filldraw[black] (\i*0.2-0.1,0) circle (1pt);
}

\foreach \i in {-5,...,-4}
{
\filldraw[green] (\i*0.1,0) circle (1pt);
}
\foreach \i in {4,...,5}
{
\filldraw[green] (\i*0.1,0) circle (1pt);
}

\foreach \i in {-4,...,5}
{
\filldraw[black] (\i*0.2-0.1,-0.9) circle (1pt);
}

\foreach \i in {-4,...,5}
{
\filldraw[black] (\i*0.2,-0.8) circle (1pt);
}

\foreach \i in {-4,...,5}
{
\filldraw[black] (\i*0.2-0.1,-0.7) circle (1pt);
}

\foreach \i in {-4,...,5}
{
\filldraw[black] (\i*0.2,-0.6) circle (1pt);
}

\foreach \i in {-4,...,5}
{
\filldraw[black] (\i*0.2-0.1,-0.5) circle (1pt);
}

\foreach \i in {-4,...,-2}
{
\filldraw[black] (\i*0.2,-0.4) circle (1pt);
}
\foreach \i in {2,...,5}
{
\filldraw[black] (\i*0.2,-0.4) circle (1pt);
}

\foreach \i in {-4,...,-3}
{
\filldraw[black] (\i*0.2-0.1,-0.3) circle (1pt);
}
\foreach \i in {4,...,5}
{
\filldraw[black] (\i*0.2-0.1,-0.3) circle (1pt);
}
\foreach \i in {-1,...,1}
{
\filldraw[black] (\i*0.2,-0.3) circle (1pt);
}

\foreach \i in {-4,...,-3}
{
\filldraw[black] (\i*0.2,-0.2) circle (1pt);
}
\foreach \i in {3,...,5}
{
\filldraw[black] (\i*0.2,-0.2) circle (1pt);
}
\foreach \i in {-1,...,2}
{
\filldraw[black] (\i*0.2-0.1,-0.2) circle (1pt);
}

\foreach \i in {-4,...,-2}
{
\filldraw[black] (\i*0.2-0.1,-0.1) circle (1pt);
}
\foreach \i in {3,...,5}
{
\filldraw[black] (\i*0.2-0.1,-0.1) circle (1pt);
}
\foreach \i in {-1,...,1}
{
\filldraw[black] (\i*0.2,-0.1) circle (1pt);
}

\end{scope}
\end{tikzpicture}
\end{subfigure}
\begin{subfigure}{0.325\textwidth}
\begin{tikzpicture}[scale=3]
\begin{scope}
\clip (-9mm,-9mm) rectangle (9mm,9 mm);
\draw[step=1mm,help lines] (-25mm,-25mm) grid (25mm,25mm);

\foreach \i in {-4,...,5}
{
\filldraw[black] (\i*0.2-0.1,0.9) circle (1pt);
}

\foreach \i in {-4,...,5}
{
\filldraw[black] (\i*0.2,0.8) circle (1pt);
}

\foreach \i in {-4,...,5}
{
\filldraw[black] (\i*0.2-0.1,0.7) circle (1pt);
}

\foreach \i in {-4,...,5}
{
\filldraw[black] (\i*0.2,0.6) circle (1pt);
}

\foreach \i in {-4,...,5}
{
\filldraw[black] (\i*0.2-0.1,0.5) circle (1pt);
}

\foreach \i in {-4,...,-2}
{
\filldraw[black] (\i*0.2,0.4) circle (1pt);
}
\foreach \i in {2,...,5}
{
\filldraw[black] (\i*0.2,0.4) circle (1pt);
}

\foreach \i in {-1,...,1}
{
\filldraw[green] (\i*0.2,0.5) circle (1pt);
}
\foreach \i in {-1,...,1}
{
\filldraw[green] (\i*0.2,-0.5) circle (1pt);
}

\foreach \i in {-5,...,-5}
{
\filldraw[green] (\i*0.1,0.4) circle (1pt);
}
\foreach \i in {-3,...,3}
{
\filldraw[green] (\i*0.1,0.4) circle (1pt);
}
\foreach \i in {5,...,5}
{
\filldraw[green] (\i*0.1,0.4) circle (1pt);
}
\foreach \i in {-5,...,-5}
{
\filldraw[green] (\i*0.1,-0.4) circle (1pt);
}
\foreach \i in {-3,...,3}
{
\filldraw[green] (\i*0.1,-0.4) circle (1pt);
}
\foreach \i in {5,...,5}
{
\filldraw[green] (\i*0.1,-0.4) circle (1pt);
}

\foreach \i in {-4,...,-3}
{
\filldraw[black] (\i*0.2-0.1,0.3) circle (1pt);
}
\foreach \i in {4,...,5}
{
\filldraw[black] (\i*0.2-0.1,0.3) circle (1pt);
}
\foreach \i in {-1,...,1}
{
\filldraw[black] (\i*0.2-0.1,0.3) circle (1pt);
}

\foreach \i in {-6,...,-4}
{
\filldraw[green] (\i*0.1,0.3) circle (1pt);
}
\foreach \i in {2,...,6}
{
\filldraw[green] (\i*0.1,0.3) circle (1pt);
}
\foreach \i in {-6,...,-4}
{
\filldraw[green] (\i*0.1,-0.3) circle (1pt);
}
\foreach \i in {2,...,6}
{
\filldraw[green] (\i*0.1,-0.3) circle (1pt);
}
\foreach \i in {-1,...,1}
{
\filldraw[green] (\i*0.2,0.3) circle (1pt);
}
\foreach \i in {-1,...,1}
{
\filldraw[green] (\i*0.2,-0.3) circle (1pt);
}

\foreach \i in {-4,...,-3}
{
\filldraw[black] (\i*0.2,0.2) circle (1pt);
}
\foreach \i in {3,...,5}
{
\filldraw[black] (\i*0.2,0.2) circle (1pt);
}
\foreach \i in {-1,...,2}
{
\filldraw[black] (\i*0.2-0.2,0.2) circle (1pt);
}

\foreach \i in {-5,...,-5}
{
\filldraw[green] (\i*0.1,0.2) circle (1pt);
}
\foreach \i in {3,...,5}
{
\filldraw[green] (\i*0.1,0.2) circle (1pt);
}
\foreach \i in {-5,...,-5}
{
\filldraw[green] (\i*0.1,-0.2) circle (1pt);
}
\foreach \i in {3,...,5}
{
\filldraw[green] (\i*0.1,-0.2) circle (1pt);
}

\foreach \i in {-4,...,-2}
{
\filldraw[black] (\i*0.2-0.1,0.1) circle (1pt);
}
\foreach \i in {3,...,5}
{
\filldraw[black] (\i*0.2-0.1,0.1) circle (1pt);
}
\foreach \i in {-1,...,1}
{
\filldraw[black] (\i*0.2-0.1,0.1) circle (1pt);
}

\foreach \i in {2,...,4}
{
\filldraw[green] (\i*0.1,0.1) circle (1pt);
}
\foreach \i in {2,...,4}
{
\filldraw[green] (\i*0.1,-0.1) circle (1pt);
}

\foreach \i in {-4,...,-3}
{
\filldraw[black] (\i*0.2,0) circle (1pt);
}
\foreach \i in {3,...,5}
{
\filldraw[black] (\i*0.2,0) circle (1pt);
}
\foreach \i in {-1,...,2}
{
\filldraw[black] (\i*0.2-0.2,0) circle (1pt);
}

\foreach \i in {3,...,5}
{
\filldraw[green] (\i*0.1,0) circle (1pt);
}

\foreach \i in {-4,...,5}
{
\filldraw[black] (\i*0.2-0.1,-0.9) circle (1pt);
}

\foreach \i in {-4,...,5}
{
\filldraw[black] (\i*0.2,-0.8) circle (1pt);
}

\foreach \i in {-4,...,5}
{
\filldraw[black] (\i*0.2-0.1,-0.7) circle (1pt);
}

\foreach \i in {-4,...,5}
{
\filldraw[black] (\i*0.2,-0.6) circle (1pt);
}

\foreach \i in {-4,...,5}
{
\filldraw[black] (\i*0.2-0.1,-0.5) circle (1pt);
}

\foreach \i in {-4,...,-2}
{
\filldraw[black] (\i*0.2,-0.4) circle (1pt);
}
\foreach \i in {2,...,5}
{
\filldraw[black] (\i*0.2,-0.4) circle (1pt);
}

\foreach \i in {-4,...,-3}
{
\filldraw[black] (\i*0.2-0.1,-0.3) circle (1pt);
}
\foreach \i in {4,...,5}
{
\filldraw[black] (\i*0.2-0.1,-0.3) circle (1pt);
}
\foreach \i in {-1,...,1}
{
\filldraw[black] (\i*0.2-0.1,-0.3) circle (1pt);
}

\foreach \i in {-4,...,-3}
{
\filldraw[black] (\i*0.2,-0.2) circle (1pt);
}
\foreach \i in {3,...,5}
{
\filldraw[black] (\i*0.2,-0.2) circle (1pt);
}
\foreach \i in {-1,...,2}
{
\filldraw[black] (\i*0.2-0.2,-0.2) circle (1pt);
}

\foreach \i in {-4,...,-2}
{
\filldraw[black] (\i*0.2-0.1,-0.1) circle (1pt);
}
\foreach \i in {3,...,5}
{
\filldraw[black] (\i*0.2-0.1,-0.1) circle (1pt);
}
\foreach \i in {-1,...,1}
{
\filldraw[black] (\i*0.2-0.1,-0.1) circle (1pt);
}

\filldraw[green] (-0.5,0) circle (1pt);
\filldraw[green] (-0.4,0.1) circle (1pt);
\filldraw[green] (-0.4,-0.1) circle (1pt);

\end{scope}
\end{tikzpicture}
\end{subfigure}
\begin{subfigure}{0.325\textwidth}
\begin{tikzpicture}[scale=3]
\begin{scope}
\clip (-9mm,-9mm) rectangle (9mm,9 mm);
\draw[step=1mm,help lines] (-25mm,-25mm) grid (25mm,25mm);

\foreach \i in {-4,...,5}
{
\filldraw[black] (\i*0.2-0.1,0.9) circle (1pt);
}

\foreach \i in {-4,...,5}
{
\filldraw[black] (\i*0.2,0.8) circle (1pt);
}

\foreach \i in {-4,...,5}
{
\filldraw[black] (\i*0.2-0.1,0.7) circle (1pt);
}

\foreach \i in {-4,...,5}
{
\filldraw[black] (\i*0.2,0.6) circle (1pt);
}

\foreach \i in {-4,...,5}
{
\filldraw[black] (\i*0.2-0.1,0.5) circle (1pt);
}

\foreach \i in {-4,...,-2}
{
\filldraw[black] (\i*0.2,0.4) circle (1pt);
}
\foreach \i in {2,...,5}
{
\filldraw[black] (\i*0.2,0.4) circle (1pt);
}

\foreach \i in {-1,...,1}
{
\filldraw[green] (\i*0.2,0.5) circle (1pt);
}
\foreach \i in {-1,...,1}
{
\filldraw[green] (\i*0.2,-0.5) circle (1pt);
}

\foreach \i in {-5,...,-5}
{
\filldraw[green] (\i*0.1,0.4) circle (1pt);
}
\foreach \i in {-1,...,2}
{
\filldraw[green] (\i*0.2-0.1,0.4) circle (1pt);
}

\foreach \i in {5,...,5}
{
\filldraw[green] (\i*0.1,0.4) circle (1pt);
}
\foreach \i in {-5,...,-5}
{
\filldraw[green] (\i*0.1,-0.4) circle (1pt);
}
\foreach \i in {-1,...,2}
{
\filldraw[green] (\i*0.2-0.1,-0.4) circle (1pt);
}
\foreach \i in {5,...,5}
{
\filldraw[green] (\i*0.1,-0.4) circle (1pt);
}

\foreach \i in {-4,...,-3}
{
\filldraw[black] (\i*0.2-0.1,0.3) circle (1pt);
}
\foreach \i in {4,...,5}
{
\filldraw[black] (\i*0.2-0.1,0.3) circle (1pt);
}
\foreach \i in {-1,...,1}
{
\filldraw[black] (\i*0.2-0.1,0.3) circle (1pt);
}

\foreach \i in {-6,...,-4}
{
\filldraw[green] (\i*0.1,0.3) circle (1pt);
}
\foreach \i in {2,...,6}
{
\filldraw[green] (\i*0.1,0.3) circle (1pt);
}
\foreach \i in {-6,...,-4}
{
\filldraw[green] (\i*0.1,-0.3) circle (1pt);
}
\foreach \i in {2,...,6}
{
\filldraw[green] (\i*0.1,-0.3) circle (1pt);
}
\foreach \i in {-1,...,1}
{
\filldraw[green] (\i*0.2,0.3) circle (1pt);
}
\foreach \i in {-1,...,1}
{
\filldraw[green] (\i*0.2,-0.3) circle (1pt);
}

\foreach \i in {-4,...,-3}
{
\filldraw[black] (\i*0.2,0.2) circle (1pt);
}
\foreach \i in {3,...,5}
{
\filldraw[black] (\i*0.2,0.2) circle (1pt);
}
\foreach \i in {-1,...,2}
{
\filldraw[black] (\i*0.2-0.2,0.2) circle (1pt);
}

\foreach \i in {-5,...,-5}
{
\filldraw[green] (\i*0.1,0.2) circle (1pt);
}
\foreach \i in {3,...,5}
{
\filldraw[green] (\i*0.1,0.2) circle (1pt);
}

\foreach \i in {-5,...,-5}
{
\filldraw[green] (\i*0.1,-0.2) circle (1pt);
}
\foreach \i in {3,...,5}
{
\filldraw[green] (\i*0.1,-0.2) circle (1pt);
}

\foreach \i in {-4,...,-2}
{
\filldraw[black] (\i*0.2-0.1,0.1) circle (1pt);
}
\foreach \i in {3,...,5}
{
\filldraw[black] (\i*0.2-0.1,0.1) circle (1pt);
}
\foreach \i in {-1,...,1}
{
\filldraw[black] (\i*0.2-0.1,0.1) circle (1pt);
}

\foreach \i in {2,...,4}
{
\filldraw[green] (\i*0.1,0.1) circle (1pt);
}

\foreach \i in {2,...,4}
{
\filldraw[green] (\i*0.1,-0.1) circle (1pt);
}

\foreach \i in {-4,...,-3}
{
\filldraw[black] (\i*0.2,0) circle (1pt);
}
\foreach \i in {3,...,5}
{
\filldraw[black] (\i*0.2,0) circle (1pt);
}
\foreach \i in {-1,...,2}
{
\filldraw[black] (\i*0.2-0.2,0) circle (1pt);
}

\foreach \i in {3,...,5}
{
\filldraw[green] (\i*0.1,0) circle (1pt);
}

\foreach \i in {-4,...,5}
{
\filldraw[black] (\i*0.2-0.1,-0.9) circle (1pt);
}

\foreach \i in {-4,...,5}
{
\filldraw[black] (\i*0.2,-0.8) circle (1pt);
}

\foreach \i in {-4,...,5}
{
\filldraw[black] (\i*0.2-0.1,-0.7) circle (1pt);
}

\foreach \i in {-4,...,5}
{
\filldraw[black] (\i*0.2,-0.6) circle (1pt);
}

\foreach \i in {-4,...,5}
{
\filldraw[black] (\i*0.2-0.1,-0.5) circle (1pt);
}

\foreach \i in {-4,...,-2}
{
\filldraw[black] (\i*0.2,-0.4) circle (1pt);
}
\foreach \i in {2,...,5}
{
\filldraw[black] (\i*0.2,-0.4) circle (1pt);
}

\foreach \i in {-4,...,-3}
{
\filldraw[black] (\i*0.2-0.1,-0.3) circle (1pt);
}
\foreach \i in {4,...,5}
{
\filldraw[black] (\i*0.2-0.1,-0.3) circle (1pt);
}
\foreach \i in {-1,...,1}
{
\filldraw[black] (\i*0.2-0.1,-0.3) circle (1pt);
}

\foreach \i in {-4,...,-3}
{
\filldraw[black] (\i*0.2,-0.2) circle (1pt);
}
\foreach \i in {3,...,5}
{
\filldraw[black] (\i*0.2,-0.2) circle (1pt);
}
\foreach \i in {-1,...,2}
{
\filldraw[black] (\i*0.2-0.2,-0.2) circle (1pt);
}

\foreach \i in {-4,...,-2}
{
\filldraw[black] (\i*0.2-0.1,-0.1) circle (1pt);
}
\foreach \i in {3,...,5}
{
\filldraw[black] (\i*0.2-0.1,-0.1) circle (1pt);
}
\foreach \i in {-1,...,1}
{
\filldraw[black] (\i*0.2-0.1,-0.1) circle (1pt);
}

\foreach \i in {2,...,3}
{
\filldraw[white] (\i*0.2-0.1,0.3) circle (1pt);
}
\foreach \i in {-1,...,1}
{
\filldraw[white] (\i*0.2,0.4) circle (1pt);
}
\foreach \i in {-1,...,1}
{
\filldraw[white] (\i*0.2,-0.4) circle (1pt);
}
\foreach \i in {2,...,3}
{
\filldraw[white] (\i*0.2-0.1,-0.3) circle (1pt);
}
\foreach \i in {2,...,2}
{
\filldraw[white] (\i*0.2,0.2) circle (1pt);
}
\foreach \i in {2,...,2}
{
\filldraw[white] (\i*0.2,-0.2) circle (1pt);
}
\foreach \i in {4,...,4}
{
\filldraw[white] (\i*0.1,0) circle (1pt);
}
\foreach \i in {3,...,3}
{
\filldraw[white] (\i*0.1,-0.1) circle (1pt);
}
\foreach \i in {3,...,3}
{
\filldraw[white] (\i*0.1,0.1) circle (1pt);
}

\foreach \i in {2,...,3}
{
\filldraw[pattern=north west lines, pattern color=blue] (\i*0.2-0.1,0.3) circle (1pt);
}
\foreach \i in {-1,...,1}
{
\filldraw[pattern=north west lines, pattern color=blue] (\i*0.2,0.4) circle (1pt);
}
\foreach \i in {-1,...,1}
{
\filldraw[pattern=north west lines, pattern color=blue] (\i*0.2,-0.4) circle (1pt);
}
\foreach \i in {2,...,3}
{
\filldraw[pattern=north west lines, pattern color=blue] (\i*0.2-0.1,-0.3) circle (1pt);
}
\foreach \i in {2,...,2}
{
\filldraw[pattern=north west lines, pattern color=blue] (\i*0.2,0.2) circle (1pt);
}
\foreach \i in {2,...,2}
{
\filldraw[pattern=north west lines, pattern color=blue] (\i*0.2,-0.2) circle (1pt);
}
\foreach \i in {4,...,4}
{
\filldraw[pattern=north west lines, pattern color=blue] (\i*0.1,0) circle (1pt);
}
\foreach \i in {3,...,3}
{
\filldraw[pattern=north west lines, pattern color=blue] (\i*0.1,-0.1) circle (1pt);
}
\foreach \i in {3,...,3}
{
\filldraw[pattern=north west lines, pattern color=blue] (\i*0.1,0.1) circle (1pt);
}
\filldraw[green] (-0.5,0) circle (1pt);
\filldraw[green] (-0.4,0.1) circle (1pt);
\filldraw[green] (-0.4,-0.1) circle (1pt);

\end{scope}
\end{tikzpicture}
\end{subfigure}
	\caption{Illustration of a configuration with one contour $\gamma$ (left in green), where the outside configuration is of type 0 and the inside configuration is of (bad) type 1. Moving the inside configuration by $-e_1$ (middle) creates a (good) configuration also inside the contour $\gamma_e$ as described in~\eqref{Con_Mov} (middle in green), however the shift also creates isolated zeros. On the right, in the large connected component of $\gamma_e$, sites are indicated in dashed blue, which can be flipped from unoccupied to occupied, as in the configuration presented in~\eqref{Ref_Con}, and therefore create an energetically preferable configuration.}
	\label{Pix_Iso}
\end{figure}
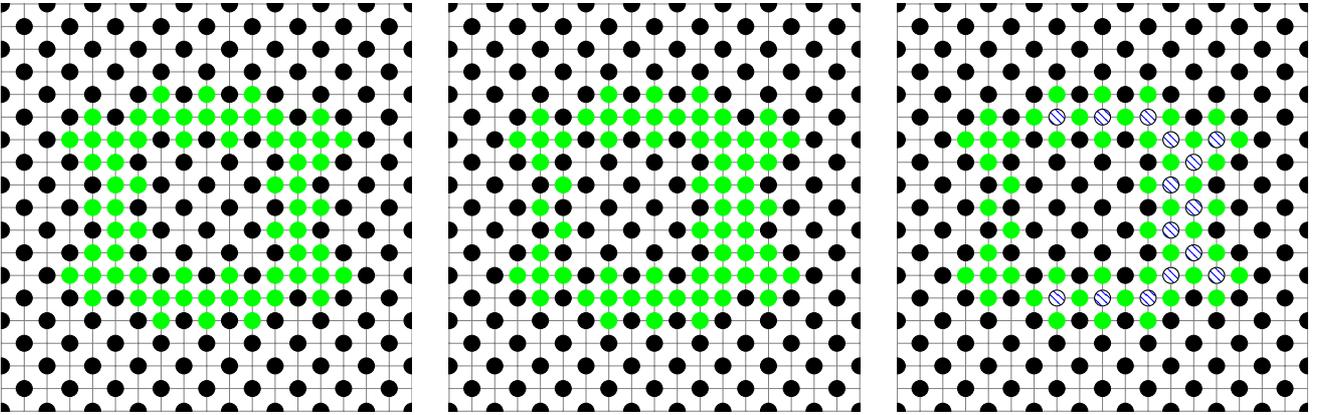
This may occur if the only neighbor of such a site is eaten up by a translated interior component. We will now explain that this is not a real problem as the number of those unwanted sites is too small to spoil our desired bound~\eqref{cd2}.

Indeed, consider the case $\Z^2$ and $e=-e_1$ for visualization. Then, in each fixed row, the sets $\gamma$ and $\gamma_e$ contain the same number of sites. To each loss site (on the left) there corresponds a site in the set $\bigcup_{k=1}^lW_l \backslash (W_l+e)$ (which is added to the contour) and which is not isolated but has a nearest neighbor in $\gamma$. From this it follows that $\gamma_e$ still has a nearest-neighbor-connected component $\tilde\gamma_e$ that is at least of size $|\gamma|/2$. Now, already $\tilde\gamma_e$ being connected, it hits a fraction (which we call $\tilde c_d>0$) of $S^0$ and we arrive at 
\begin{equation}\label{cd3}
\begin{split}
&|\gamma_e \cap S^0|\geq |\tilde \gamma_e \cap S^0|
\geq \tilde c_d |\tilde \gamma_e| \geq c_d|\gamma|, 
\end{split}
\end{equation}
where $c_d=\tilde c_d/2$. Note that we may regard $\tilde\gamma_e$, as obtained in this process, as a contour, and its removal in the line of Case 1. 
 
This proves the claim in the general form and hence proves the Peierls estimate. 
\end{proof}

\subsubsection{Non-removable discontinuities for $\mu_p'$ via the relation to Bernoulli fields conditioned on isolation}\label{Sec_Non}
\begin{proof}[Proof of Theorem~\ref{thm_large_p}] Suppose that $\gamma'$ is a any specification for $\mu_p'$. 
Consider a square $Q$ of sidelength $3$ around the origin. It will play the role of an observation window. Recall that we write $\omega'$ for the second-layer variables, and we write $\omega$ for the first-layer variables. 
We will draw the assumption that 
$$\omega'_{Q^c}\mapsto \gamma'_Q(\d \omega'_Q|\omega'_{Q^c})$$ 
is continuous at $\omega'_{Q^c}=0'_{Q^c}$ to a contradiction, by exhibiting a finite jump size. 
We use the notation $1'_{B}$ for the all-one configuration in the volume $B$. A single-site observation window would not show the phenomenon, as a conditioning $0'_{0^c}$ for the model conditional on non-isolation forces the origin to be $0'$ due to the non-isolation constraint.

For the purpose of showing the persistence of jumps on $Q$, we consider the loophole volumes $B_L$ and $B_L+e$ as above, and choose for each $L$ cubes $C_L$ that contain both $B_L$ and $B_L+e$ with a layer of sufficiently large finite thickness. Thickness two will do.  

In the first step, we 
consider conditional probabilities of the second-layer measure $\mu_p'=T\mu_p$ of the form
\begin{equation*}
\begin{split}
&\mu_p'(\omega'_Q| \omega'_{C_L\backslash Q}). 
\end{split}
\end{equation*}
We want to show that they are essentially discontinuous at the empty conditioning, so they cannot come from a quasilocal specification $\gamma'$. The latter argument will be discussed below in the second step, let us now discuss how to obtain the essential discontinuity. There is a slightly tricky part, as we need to go to higher volumes than single-sites, and need to take care of the constraints very carefully. 
In order to do consider
\begin{equation*}
\begin{split}
&\mu_p'(\omega'_Q| 0'_{B_L\backslash Q}1'_{C_L\backslash B_L}) 
\end{split}
\end{equation*}
where $\omega'_Q\in \{0,1\}^Q$. 

It is useful to compare with the empty configuration on the second layer
\begin{equation*}
\begin{split}
&\frac{\mu_p'(\omega'_Q| 0'_{B_L\backslash Q}1'_{C_L\backslash B_L})}
{\mu_p'(0'_Q| 0'_{B_L\backslash Q}1'_{C_L\backslash B_L})}=\frac{\mu_p'(\omega'_Q 0'_{B_L\backslash Q}1'_{C_L\backslash B_L})}
{\mu_p'(0'_Q 0'_{B_L\backslash Q}1'_{C_L\backslash B_L})}.
 \end{split}
\end{equation*}
The denominators never vanish, as the configuration that appears in the denominator on the right-hand side obeys the non-isolation constraint. 

Given the boundary condition $0'_{B_L\backslash Q}$, the non-isolation constraint on the second layer limits the configurations to non-isolated configurations on the cube $Q$ in order to have a non-zero measure. 

We will take as a useful local pattern on $Q$, the following second-layer reference configuration, which is most easily visualized in two dimensions where it looks as follows 
 $$\omega'^*_Q=
 \begin{array}{rrr} 
0 & 1 & 0 \\ 
1 & 1 & 1 \\ 
0 & 1 & 0 \\ 
\end{array}.$$
It is clearly compatible with the second-layer constraint as $\omega'^*_Q 0'_{Q^c}$ contains 
no isolated ones.  
In general dimensions $d$ we choose it analogously,  namely as the  checkerboard groundstate of type-0 but with an additional $1$ at the origin, i.e., 
\begin{equation}\label{Eq_Om}
\begin{split}
(\omega'^*_Q)_i=
\begin{cases}
\omega^0_i & i \in Q\backslash 0 \\ 
1 & i=0.
\end{cases}
 \end{split}
\end{equation}
We have the following useful observation. Given $\omega'^*_Q$, the underlying Bernoulli-field configuration $\omega_Q$ from which it appears as $T$-image, must take the same values on $Q$, independently of $\omega'_{Q^c}$.  

To understand the following steps it will be helpful to make a notational distinction and write $\omega^*_Q=\omega'^*_Q$ when we refer to the same configuration as a {\em first-layer} configuration. With this we have 
\begin{equation*}
\begin{split}
&\frac{\mu_p'(\omega'^*_Q 0'_{B_L\backslash Q}1'_{C_L\backslash B_L})}
{\mu_p'(0'_Q 0'_{B_L\backslash Q}1'_{C_L\backslash B_L})}=\frac{\mu_{p,C_L}(T\s_{C_L}=\omega'^*_Q 0'_{B_L\backslash Q}1'_{C_L\backslash B_L})}
{\mu_{p,C_L}(T\s_{C_L}=0'_Q 0'_{B_L\backslash Q}1'_{C_L\backslash B_L})}.
 \end{split}
\end{equation*}
Next, we are aiming for a reformulation that involves only quantities of the first-layer model with isolation constraint in the {\em whole volume} $B_L$ including $Q$. For this we perform some manipulations. 
We split the numerator as follows,
\begin{equation*}
\begin{split}
&\mu_{p,C_L}(T\s_{C_L}=\omega'^*_Q 0'_{B_L\backslash Q}1'_{C_L\backslash B_L})=\mu_{p,Q}(\omega^*_Q) \mu_{p,C_L \backslash Q}\Bigl(T(\omega^*_Q,\s_{C_L\backslash Q})\big|_{C_L\backslash Q} =0'_{B_L\backslash Q}1'_{C_L\backslash B_L}\Bigr).
 \end{split}
\end{equation*}
Now, we change  the middle site on $Q$ on the right-hand side from $1$ to $0$ to obtain a first-layer configuration that obeys the isolation constraint on $Q$. For this, we first write the simple identity
$$\mu_{p,C_L}(\omega^*_Q)=\tfrac{p}{1-p}
\mu_{p,C_L}(\omega^0_Q).$$ 
It is important to note that we may also replace 
$$T_{C_L}(\omega^*_Q,\s_{C_L\backslash Q})\big|_{C_L\backslash Q}
=T_{C_L}(\omega^0_Q,\s_{C_L\backslash Q})\big|_{C_L\backslash Q},$$
which is possible as the middle site of $Q$ has no influence on the values of the restriction of 
$T_{C_L}(\omega^0_Q,\s_{C_L\backslash Q})$ to $Q^c$. 
So we arrive 
at 
\begin{equation*}
\begin{split}
&\mu_{p,C_L}(T\s_{C_L}=\omega'^*_Q 0'_{B_L\backslash Q}1'_{C_L\backslash B_L})=\frac{p}{1-p}
\mu_{p,C_L}(\s_Q=\omega^0_Q, T_{C_L}(\omega^0_Q,\s_{C_L\backslash Q})\big|_{C_L\backslash Q}
=0'_{B_L\backslash Q}1'_{C_L\backslash B_L}).
 \end{split}
\end{equation*}
Now we have achieved our goal, as we may recognize that
\begin{equation*}
\begin{split}
&\frac{\mu_{p,C_L}(\s_Q=\omega^0_Q, 
T_{C_L}(\omega^0_Q,\s_{C_L\backslash Q})\big|_{C_L\backslash Q}
=0'_{B_L\backslash Q}1'_{C_L\backslash B_L})}
{\mu_{p,C_L}(T\s_{C_L}=0'_Q 0'_{B_L\backslash Q}1'_{C_L\backslash B_L})}=\nu_{B_L}(\sigma_Q=\omega^0_Q),
 \end{split}
\end{equation*}
with the conditional first-layer measure $\nu_{B_L}$ that is conditioned on isolation, in the whole volume $B_L$, as defined in the beginning of Section~\ref{Sec_Pei}. Indeed, the denominator on the l.h.s.~is the partition function of the conditional measure, up to the terms on $C_L\backslash B_L$ on which the first-layer configuration is frozen, and which cancel against those in the numerator. In particular there is no $C_L$-dependence. 
We have thus proved the following {\em unfixing lemma.}
\begin{lem}\label{lem_NonGibbs}
The second-layer conditional probabilities and the first-layer model under the non-isolation constraint satisfy the following relation 
\begin{equation*}
\begin{split}
&\frac{\mu_p'(\sigma_Q'=\omega'^*_Q| 0'_{B_L\backslash Q}1'_{C_L\backslash B_L})}
{\mu_p'(\sigma_Q'=0'_Q| 0'_{B_L\backslash Q}1'_{C_L\backslash B_L})}=\frac{p}{1-p}\nu_{B_L}(\sigma_Q=\omega^0_Q),
 \end{split}
\end{equation*}
where $\omega^0_Q$ denotes the checkerboard groundstate and $\omega'^*_Q$ is defined in~\eqref{Eq_Om}. 
The same relation holds for the shifted volume $B_L+e$.  
\end{lem}
Note that this representation is particularly nice, as we are reduced to the discussion of the first-layer model in the full volume $B_L$ (and not a volume reduced by $Q$) where we have the Peierls estimate to our disposition. 
In particular, by the Peierls estimate we have that
$$\nu_{B_L}(\sigma_Q=\omega^0_Q)\geq 1-\sum_{j\in Q}
\mu_{B_L}(\sigma_j\neq \omega^0_j)\geq 1-|Q|\epsilon(p),$$
while 
$$\nu_{B_L+e}(\s_Q=\omega^0_Q)\leq \epsilon(p).$$

In the final step, we bring the arbitrarily chosen specification $\gamma'$ into play, with the aim to show that it must inherit a discontinuity at the empty configuration, too. For any pattern $\o'_\Q$ in the observation window $Q$, we bound the infimum over perturbations of the empty configuration outside the volume $\Delta_L:=B_L\cap B_{L+e}$ via 
\begin{equation*}
\begin{split}
\mu_p'(\omega'_Q| 0'_{B_L\backslash Q}1'_{C_L\backslash B_L})&=\int \gamma'_Q(\omega'_Q| 0'_{B_L\backslash Q}1'_{C_L\backslash B_L}
\tilde\omega'_{C_L^c})\,\,\mu_p'(\d\tilde\omega'_{C_L^c}| \omega'_{C_L\backslash Q}=0'_{B_L\backslash Q}1'_{C_L\backslash B_L})\cr
&\geq \inf_{\omega'_{\Delta_L^c}}
\gamma'_Q(\omega'_Q| 0'_{\Delta_L\backslash Q}
\omega'_{\Delta_L^c})=:a_L(\omega'_Q).
\end{split}
\end{equation*}
Similar arguments give that
\begin{equation*}
\begin{split}
\mu_p'(\omega'_Q| 0'_{B_L+e\backslash Q}1'_{C_L \backslash B_L+e}) &\geq a_L(\omega'_Q)\cr
\mu_p'(\omega'_Q| 0'_{B_L\backslash Q}1'_{C_L\backslash B_L}) &\leq \sup_{\omega'_{\Delta_L^c}}
\gamma_Q(\omega'_Q| 0'_{\Delta_L\backslash Q}
\omega'_{\Delta_L^c})=:b_L(\omega'_Q)\cr
\mu_p'(\omega'_Q| 0'_{B_L+e\backslash Q}1'_{C_L\backslash B_L+e})&\leq b_L(\omega'_Q).
\end{split}
\end{equation*}
Now consider specifically the patterns $0'_Q$ and $\omega'^*_Q$ and note that
\begin{equation}
\begin{split}\label{uppereast}
&\frac{p}{1-p}\nu_{B_L}(\sigma_Q=\omega^0_Q)=\frac{\mu_p'(\omega'^*_Q| 0'_{B_L\backslash Q}1'_{C_L\backslash B_L})}
{\mu_p'(0'_Q| 0'_{B_L\backslash Q}1'_{C_L\backslash B_L})}\leq \frac{b_L(\omega'^*_Q)}{a_L(0_Q)}
 \end{split}
\end{equation}
and 
\begin{equation}
\begin{split}\label{lowereast}
&\frac{p}{1-p}\nu_{B_L+e}(\sigma_Q=\omega^0_Q)=\frac{\mu_p'(\omega'^*_Q| 0'_{B_L+e\backslash Q}1'_{C_L\backslash B_L+e})}
{\mu_p'(0'_Q| 0'_{B_L+e\backslash Q}1'_{C_L\backslash B_L+e})}
\geq \frac{a_L(\omega'^*_Q)}{b_L(0_Q)},
 \end{split}
\end{equation}
and remark that the denominators are uniformly bounded against zero. 

Now, by the Peierls estimate presented in Lemma~\ref{lem_Pei} in previous section, we have lower bounds on the left-hand side of~\eqref{uppereast} and upper bounds on the left-hand side of~\eqref{lowereast}. These contradict the continuity assumption on the specification $\gamma'$, i.e., that the right-hand sides have the same limit as $L\uparrow\infty$. 

This proves the discontinuity of the specification kernel $\gamma'_Q$ for any arbitrary specification $\gamma'$, at the fully empty configuration. 
\end{proof}

\subsection{Proof of Theorem~\ref{thm_small_p}: Gibbsianness}
\label{Sec_Proofs_Low_p_s}
In this section, we construct a continuous specification $\g'$ for $\mu_p'$ for small $p$. The main ingredient is an application of the Dobrushin-uniqueness bound and the backward-martingale theorem. In the first step, we construct the conditional probabilities in finite volumes. 

For this we use the following notation. For $\L\subset\Z^d$ we denote by $\L^c:=\Z^d\setminus \L$ its {\em complement} and by $\partial_-\L:=\{x\in \L\colon \text{there exists }y\in \L^c\text{ with }y\sim_{\Z^d} x\}$ its {\em interior boundary}. The set $\L^o:=\L\setminus\partial_-\L$ denotes the {\em interior} and $\bar\L:=((\L^c)^o)^c$ the {\em extension} of $\L$. Moreover, $\partial_+\L:=\bar\L\setminus \L$ then denotes the {\em outer boundary} of $\L$.

\subsubsection{The specification}
Let us consider, for any (large) finite volume $\D\Subset\Z^d$, conditional probabilities of $\mu'_p$ inside $\D$ also given a (first-layer) boundary condition $\o$ outside $\D$. More precisely, let $\L\subset\D$, $\o'=\o'_\L\o'_{\D\setminus\L}\o'_{\D^c}\in \O'$ and let $\o$ be such that $\o_{\D^c}\in T^{-1}(\o'_{(\D^o)^c})$, then 
\begin{equation}\label{Eq0}
\begin{split}
\g'_{\o_{\D^c}, \L}&(\o'_\L|\o'_{\D\setminus \L})
:=\frac{
\sum_{\tilde\o_{\D}}\mu_p(\tilde\o_{\D})\one\{T_\D(\tilde\o_{\D} \o_{\D^c})=\o'_\D\}}{
\sum_{\tilde\o_{\D\setminus\L^o}}\mu_p(\tilde\o_{\D\setminus\L^o})\one\{T_{\D \setminus \L}(\tilde\o_{\D\setminus\L^o} \o_{\D^c})=\o'_{\D\setminus\L}\}
}\\
=&\frac{
\sum_{\tilde\o_{\D\setminus\L^o}}\mu_p(\tilde\o_{\D\setminus\L^o})\one\{T_{\D \setminus \L}(\tilde\o_{\D\setminus\L^o} \o_{\D^c})=\o'_{\D\setminus\L}\}
\sum_{\tilde\o_{\L^o}}\mu_p(\tilde\o_{\L^o})\one\{T_\L(\tilde\o_{\L^o}\tilde\o_{\D \setminus \L^o} )=\o'_{ \L} \}}{
\sum_{\tilde\o_{\D\setminus\L^o}}\mu_p(\tilde\o_{\D\setminus\L^o})\one\{T_{\D \setminus \L}(\tilde\o_{\D\setminus\L^o} \o_{\D^c})=\o'_{\D\setminus\L}\}
}\\
=&\frac{
\sum_{\tilde\o_{\D\setminus\L^o}}\mu_p(\tilde\o_{\D\setminus\L^o})\one\{T_{\D \setminus \L}(\tilde\o_{\D\setminus\L^o} \o_{\D^c})=\o'_{\D\setminus\L}\}
f_{\o'_\L}(\tilde\o_{\partial_- \L\cup\partial_+ \L }) }{
\sum_{\tilde\o_{\D\setminus\L^o}}\mu_p(\tilde\o_{\D\setminus\L^o})\one\{T_{\D \setminus \L}(\tilde\o_{\D\setminus\L^o} \o_{\D^c})=\o'_{\D\setminus\L}\}
},
\end{split}
\end{equation}
where we wrote $T_\L(\o)$ instead of $(T(\o))_\L$ and 
$$
f_{\o'_\L}(\o_{\partial_- \L\cup\partial_+ \L }):= \sum_{\tilde\o_{\L^o}}\mu_p(\tilde\o_{\L^o})\one\{T_{\L}(\tilde\o_{\L^o}\o_{\partial_- \L\cup\partial_+ \L })=\o'_{ \L} \}
$$
is a local function. We have the following consistency result. 
\begin{lem}\label{lem_specification}
Assume that, given $\L \Subset \Z^d$ and $\o' \in \O'$, $\lim_{\D \uparrow \zd} \g'_{\o_{\D^c}, \D}(\o'_\L|\o'_{\D\sm \L})=:\g'_\L(\o'_\L | \o'_{\L^c})$ exists and is independent of $\o_{\D^c} \in T^{-1}(\o'_{(\D^o)^c})$. Then, $\g'$ is a specification for $\mu'_p$.
\end{lem}
\begin{proof}
First note that for any $\o'\in \O'$ we can estimate, 
\begin{equation}\label{Eq_o}
\begin{split}
\bigl\lvert \mu'_p(\o'_\L | \o'_{\D\sm \L}) -  \g'_\L(\o'_\L | \o'_{\L^c})\bigr \rvert
&\leq \sup_{\o_{\partial_+\D}} \bigl \lvert \g'_{\o_{\partial_+\D}, \D}(\o'_\L | \o'_{\D\sm \L})  -  \g'_\L(\o'_\L | \o'_{\L^c})\bigr \rvert,
\end{split}
\end{equation}
where the supremum is taken over suitable boundary configurations compatible with $\o'$, since $\mu'_p(\o'_\L | \o'_{\D\sm \L})$ can be written as an integral with respect to $\g'_{\cdot, \D}(\o'_\L | \o'_{\D\sm \L})$. In particular, under our assumptions, the right-hand side of Equation~\ref{Eq_o} tends to zero as $\D$ tends to $\Z^d$. 

Now, consider a cofinal sequence $\D_n \uparrow \L^c$ and let $(\FF'_{\D_n})_{n \in \N}$ denote the canonical filtration on $\O'$. Note that the sequence of random variables $\mu'_p(\o'_\L | \FF'_{\D_n})$ that are $\mu'_p$-almost surely defined as $\mu'_p(\o'_\L | \FF'_{\D_n})(\o') =  \mu'_p(\o'_\L | \o'_{\D_n \setminus \L})$, is a uniformly-integrable martingale adapted to $(\FF'_{\D_n})_{n \in \N}$. Since $\s \left( \bigcup_{n} \FF'_{\D_n} \right) = \FF'_{\L^c}$, by L\'evy's zero-one law, $(\mu'_p(\o'_\L | \FF'_{\D_n}))_{n \in \N}$ converges $\mu'_p$-almost surely and in $L^1$ towards $\mu'_p(\o'_\L | \FF'_{\L^c})$ as $n$ tends to infinity. But this implies that for any $\tilde\o'_\L$, we can pick $n$ sufficiently large such that 
\begin{equation*}
\begin{split}
&\int\mu'_p(\o')|\mu'_p(\tilde\o'_\L|\o'_{\L^c})-\g'_\L(\tilde\o'_\L|\o'_{\L ^c})|\\
&\le \int\mu'_p(\o')|\mu'_p(\tilde\o'_\L|\o'_{\L^c})-\mu'_p(\tilde\o'_\L|\o'_{\D_n\sm\L})|+ \int\mu'_p(\o')|\mu'_p(\tilde\o'_\L|\o'_{\D_n\sm\L})-\g'_\L(\tilde\o'_\L|\o'_{\L ^c})|<\eps,
\end{split}
\end{equation*}
where we used L\'evy's zero-one law in the first summand and the bound~\eqref{Eq_o} in the second summand on the right-hand side. But this implies that $\int\mu'_p(\o'_{\L^c})\g'_\L(\tilde\o'_\L|\o'_{\L ^c})=\mu'_p(\tilde\o_\L')$ and hence $\g'$ is a specification for $\mu'_p$. 
\end{proof}

\subsubsection{Transformations into first-layer constraint models}
\begin{figure}[!htpb]
\centering
	\input{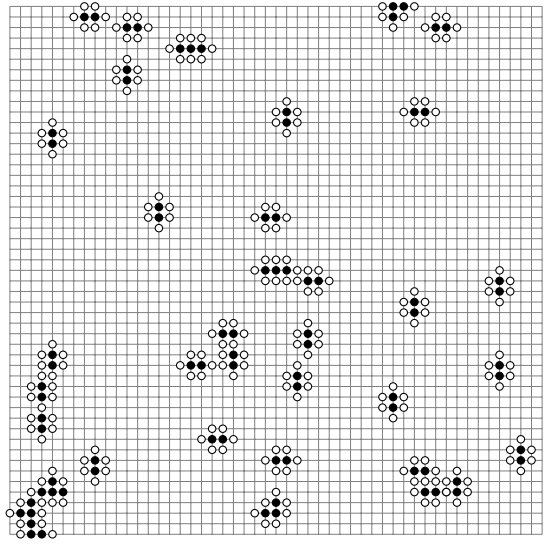}
	\caption{Illustration of the fixed area (black and white dots) based on a thinned configuration (black dots). The thinned configuration is surrounded by unoccupied sites (white dots).}
	\label{Pix_Fix}
\end{figure}

In order to establish the conditions of Lemma~\ref{lem_specification} for sufficiently small $p$, we employ the Dobrushin uniqueness theorem for the first-layer constraint model as defined in~\eqref{eq_2nd layer1}. 
For this, first note that we can uniquely identify $\o'$ with the subset of its occupied sites in $\Z^d$ and with some notational abuse $\bar\o'\subset\Z^d$ of $\o'$ is a {\em fixed area} in the sense that, under $T$, there is no choice for the Bernoulli field in how to realize $\o'$. Recall that $\o'$ consist of clusters of size at least two and $\bar\o'$ then consists of clusters of size at least two surrounded by unoccupied sites, see Figure~\ref{Pix_Fix} for an illustration. 

In view of this, we introduce the following specification associated to the {\em first-layer constraint model} on $\O$
\begin{equation}\label{eq_2nd layer1}
\g^{S}_{\D}(\o_\D|\o_{\D^c})
:=\frac{
\mu_p(\o_{\D\cap S})\one\{\o_{\D\cap S}\o_{\D^c\cap S}\text{ is $T$-feasible on } \D\cap S\}}{
\sum_{\ot_{\D\cap S}}\mu_p(\ot_{\D\cap S})\one\{\ot_{\D\cap S}\o_{\D^c\cap S}\text{ is $T$-feasible on } \D\cap S \}}. 
\end{equation}
Here, $S\subset\Z^d$ is an {\em unfixed area} that is arbitrary at this stage, $\D\Subset \Z^d$ and any configuration $\o\in \O$ is called {\em $T$-feasible} on a set $\D\cap S$ if all occupied sites of $\o$ in $\D\cap S$ have no neighboring occupied sites in $\bar\D\cap S$. In particular, with this definition, 
\begin{equation*}
\g'_{\o_{\D^c}, \L}(\o'_\L|\o'_{\D\setminus \L})=\g^{S}_{\D}(f_{\o'_\L}|\o_{\D^c}) 
\end{equation*}
for the particular choice of the unfixed area given by 
$S=S(\o'_{\D\setminus \L})=(\D\setminus \L^o)\setminus \bar\o_{\D\setminus \L}'$. 
Here we used that, in the fixed area $\bar\o_{\D\setminus \L}'$, the Bernoulli field is completely determined by $\o_{\D\setminus \L}'$ and hence the corresponding factor cancels in~\eqref{Eq0}. The following result verifies the conditions of Lemma~\ref{lem_specification} for sufficiently small $p$.
\begin{lem}\label{lem_Dob}
Let $p<1/(2d)$. Then, for any $\L \Subset \Z^d$ and $\o' \in \O'$, the limit $\lim_{\D \uparrow \zd} \g^{S(\o'_{\D\setminus\L})}_\D(f_{\o'_\L}|\o_{\D^c})$ exists and is independent of $\o_{\D^c} \in T^{-1}(\o'_{(\D^o)^c})$.
\end{lem}
\begin{proof}
We use the Dobrushin-uniqueness approach for the specification $\g_\D^S$ as defined in~\eqref{eq_2nd layer1}, where $S\subset\Z^d$ is any unfixed area. 
Consider the Dobrushin matrix 
$$
C_{ij}(p) = \max_{\o_{j^c} = \ot_{j^c}}\|\g^S_i(\cdot |\o_{i^c})
-\g^S_i(\cdot |\ot_{i^c})\|_{\text{TV}}
$$ for $i,j\in S$, where complements are defined in $S$. Note that the exterior boundary $\partial_+S$ of $S$ consists of unoccupied sites, see Figure~\ref{Pix_Fix}. We have $C_{ij}(p)=0$ unless $i$ and $j$ are neighbors in $S$. Otherwise,
\begin{align*}
C_{ij}(p)=\tfrac{1}{2}\max_{\o_{j^c} = \ot_{j^c}}(|\g^S_i(0|\o_{i^c})
-\g^S_i(0 |\ot_{i^c})|+ |\g^S_i(1|\o_{i^c})
-\g^S_i(1 |\ot_{i^c})|)=p,
\end{align*}
where the maximum is realized when $\o_{j^c}$ is unoccupied while $\o_{j} $ is unoccupied and $\ot_{j}$ is occupied. 
In particular, for the Dobrushin criterion, we have 
$$
c(p)=\sup_{i\in S}\sum_{j \sim i}C_{ij}(p)\le 2dp,
$$
independent of $S$. By~\cite[Theorem 8.7]{Ge11}, for all $p< 1/(2d)$  and $S$, $\g_\D^S$ admits a unique infinite-volume Gibbs measure $\mu^S$. Finally, using the remark made above~\cite[Equation 8.25]{Ge11}, $\g^{S}_\D(f_{\o'_\L}|\o_{\D^c})$ converges uniformly in $\o$ towards $\mu^{S}(f_{\o'_\L})$, which finishes the proof.
\end{proof}

\subsubsection{Quasilocality of the specification}
What remains to be done in order to finish the proof of Theorem~\ref{thm_small_p} is to establish quasilocality for the specification. Let $s$ denote the $\ell_\infty$ metric on $\Z^d$ and define $s(\L,\D)=\inf\{s(i,j)\colon i\in \L, j\in \D\}$.
\begin{lem}\label{lem_Quasi}
For $p<1/(2d)$ there exist constants $C,c>0$ such that for all $\L\subset\D\Subset\Z^d$ and all configurations $\o'$ and $\tilde\o'$ with $\o'_\D=\tilde\o'_\D$ we have that 
\begin{equation*}
|\g'_\L(\o'_\L|\o'_{\L^c})-\g'_\L(\o'_\L|\tilde\o'_{\L^c})|\leq  C |\Lambda| e^{-c s(\Lambda, \D^c)}.
\end{equation*}
In particular, the specification $\g'$ is quasilocal. 
\end{lem}
\begin{proof}
We use the representation of $\g'_\L(\o'_\L|\o'_{\L^c})$ in terms of the unique infinite-volume Gibbs measure $\mu^{S}(f_{\o'_\L})$ as presented in the proof of Lemma~\ref{lem_Dob}. This representation exists since we work in the Dobrushin-uniqueness regime. Now, for the quasilocality, we use the criterion~\cite[Remark 8.26]{Ge11} applied to~\cite[Theorem 8.20]{Ge11}. More precisely, since $p<1/(2d)$, by~\cite[Theorem 8.20]{Ge11}, for $S\cap\D=S'\cap \D$, we have that 
\begin{equation*}
|\mu^{S}(f_{\o'_\L})-\mu^{S'}(f_{\o'_\L})|\le D(\L,\D), 
\end{equation*}
where $D(\L,\D)=\sum_{i\in \L, j\in \D^c}\big(\sum_{n\ge 0}C^n\big)_{i,j}$ with $C^n=C^n(p)$ the $n$-th power of the Dobrushin matrix as presented in the proof of Lemma~\ref{lem_Dob}.  Now choose $c>0$ sufficiently small such that $p \e^c<1/(2d)$, then, by~\cite[Remark 8.26]{Ge11}, 
\begin{equation*}
D(\L,\D)\le |\L|(1-2dp\e^c)^{-1}\e^{-c d(\L,\D^c)}. 
\end{equation*} 
This finishes the proof. 
\end{proof}

\section*{Acknowledgements}
We thank Nils Engler for inspiring discussions. This work was funded by the German Research Foundation under Germany's Excellence Strategy MATH+: The Berlin Mathematics Research Center, EXC-2046/1 project ID: 390685689 and the German Leibniz Association via the Leibniz Competition 2020.

\begin{scriptsize}
\bibliography{../Jahnel}
\bibliographystyle{alpha}
\end{scriptsize}
\end{document}